\documentclass[11pt,a4paper,reqno]{amsart}

\usepackage{mathabx}

\usepackage[dvipsnames]{xcolor}
\usepackage{tikz}
\usepackage{amsmath,bm,bbm,amsthm, amssymb}
\usepackage{fullpage}
\usepackage{color}
\usepackage{dsfont}
\usepackage{animate}
\usepackage{etoolbox}
\usepackage{ulem}
\usepackage{comment}
\usepackage{pgf}
\usetikzlibrary{calc}
\usetikzlibrary{patterns}
\usetikzlibrary{arrows}
\usetikzlibrary{decorations.pathreplacing}
\usepackage[utf8]{inputenc}
\usepackage{pgfplots}
\usepackage{tcolorbox}
\usepackage{adjustbox}
\usepackage[final]{pdfpages}
\usepackage[ruled,vlined,linesnumbered,algosection,resetcount]{algorithm2e}
\usepackage{blkarray}

\usepackage{diagbox}
\def\E{\mathbb E}

\def\mc{\mathcal}

\def\P{\mathbb P}
\def\R{\mathbb R}
\def\wc{\widecheck}
\def\Rm{\wc{\mathbb R}^p}
\def\xm{\wc x}
\def\cx{\xm}
\def\cz{\wc z}
\def\Bm{\wc B}
\def\Wm{\wc W}
\def\Qm{\wc Q}
\def\Q{\mathbb Q}
\def\tf{\tfrac}
\def\XX{\mc X}

\def\cX{\wc X}
\def\Z{\mathbb Z}
\def\bel{\begin{lemma}}
\def\enl{\end{lemma}}
\def\been{\begin{enumerate}}
\def\bee{\begin{example}}
\def\bepr{\begin{proposition}}
\def\bep{\begin{proof}}
\def\bet{\begin{theorem}}
\def\bede{\begin{definition}}
\def\ende{\end{definition}}

\def\bi{\big}
\def\bs{\boldsymbol}
\def\a{\alpha}
\renewcommand\b{\beta}

\def\co{\colon}

\def\d{{\rm d}}
\def\enen{\end{enumerate}}
\def\ene{\end{example}}
\def\endo{\end{document}}
\def\enpr{\end{proposition}}
\def\enp{\end{proof}}
\def\ent{\end{theorem}}

\def\es{\varnothing}
\def\e{\varepsilon}
\def\ff{\infty}
\def\f{\frac}
\def\g{\gamma}
\def\im{\item}

\def\ms{\mathsf}
\def\ni{\noindent}
\def\one{\mathbbmss{1}}

\def\r{\rho}
\def\su{\subseteq}
\def\s{\sigma}

\def\ti{\times}
\def\t{\tau}

\def\tbm{b_-}
\def\tbp{b_+}
\def\tdm{d_-}
\def\tdp{d_+}
\def\trp{r_+}
\def\trm{r_-}
\def\rfc{r_{\ms I}}
\def\rci{r_{\ms C}}
\def\rce{r_{\ms{Cech}}}

\def\sm{\setminus}

\def\ba{\,|\,}
\def\dist{\ms{dist}}
\def\ot{[0, S]^2}

\def\cv{C_{\ms V}}
\def\ca{C_{\ms A}}

\def\BB{\mc B}
\def\vp{\varphi}
\def\PP{\mc P}
\def\RAi{R^{\ms {A}}}
\def\Rci{R^{ \ms{int}}}
\def\Rce{R^{\ms{ext}}}

\def\YY{\mc Y}

\def\Var{\ms{Var}}

\def\De{\Delta}
\def\Cov{\ms{Cov}}
\def\wt{\widetilde}
\def\sb{\s_{\ms B}}
\def\sd{\s_{\ms D}}
\def\bef{\begin{figure}}
\def\enf{\end{figure}}
\def\xx{\bs x}
\def\A{\mathbb A}
\def\GG{\mc G}
\def\on{\{0, \dots, n - 1\}^p}
\def\diam{\ms{diam}}
\def\AA{\mc A}

\def\enit{\end{itemize}}
\def\TA{T_{\ms{Area}}}
\def\TF{T_{\ms I}}
\def\TP{T_{\ms{Pers}}}

\def\bbnd{\b_{n, p - 2}^{M, b, d}}
\def\bbd{\b_{p - 2}^{M, b, d}}

\theoremstyle{plain}
\newtheorem{theorem}{Theorem}%[chapter]
\newtheorem{proposition}[theorem]{Proposition}

\newtheorem{lemma}[theorem]{Lemma}

\theoremstyle{definition}
\newtheorem{definition}[theorem]{Definition}

\newtheorem{example}[theorem]{Example}

\theoremstyle{remark}

\newtheorem*{definition*}{{\bf Definition}}


\begin{document}

%\title{Goodness-of-fit tests through persistence for tessellations based on topological data analysis }
\title{Persistent homology based goodness-of-fit tests for spatial tessellations}

\author{Christian Hirsch}
\address[Christian Hirsch]{Department of Mathematics, Aarhus University, Ny Munkegade 118, 8000 Aarhus C, Denmark}
\address[Christian Hirsch]{DIGIT Center, Aarhus University, Finlandsgade 22, 8200 Aarhus N, Denmark}
\email{hirsch@math.au.dk}
\author{Johannes Krebs}
\address[Johannes Krebs]{Department of Mathematics, KU Eichst\"att-Ingolstadt, Ostenstra\ss e 28, 85072 Eichst\"att, Germany} 
\email{johannes.krebs@ku.de}
\author{Claudia Redenbach  }
\address[Claudia Redenbach]{Department of Mathematics, University of Kaiserslautern, Gottlieb-Daimler-Stra\ss e 47, 67663 Kaiserslautern,
Germany}
\email{redenbach@mathematik.uni-kl.de}

\begin{abstract}
	Motivated by the rapidly increasing relevance of virtual material design in the domain of materials science, it has become essential to assess whether topological properties of stochastic models for a spatial tessellation are in accordance with a given dataset. Recently, tools from topological data analysis such as the persistence diagram have allowed to reach profound insights in a variety of application contexts. In this work, we establish the asymptotic normality of a variety of test statistics derived from a tessellation-adapted refinement of the persistence diagram. Since in applications,  it is common to work with tessellation data subject to interactions, we establish our main results for Voronoi and Laguerre tessellations whose generators form a Gibbs point process. We elucidate how these conceptual results can be used to derive goodness of fit tests, and then investigate their power in a simulation study. Finally, we apply our testing methodology to a tessellation describing real foam data.
\end{abstract}

\vspace{4cm}
\maketitle
\section{Introduction}
\label{intr_sec}

%
%VIRTUAL MATERIAL DESIGN
%
Accurate stochastic-geometry models for the microstructure of complex materials have become an indispensable tool in modern materials science, as they open the door towards virtual material design. More precisely, by being able to rapidly generate plausible realizations of microstructures stemming from a large variety of production parameters, it is possible to reduce costly experiments to the most promising candidates for a desired functional property of the new material \cite{NEUMANN2020211,stiff,WESTHOFF20181}.

%
%TEST
%
Clearly, the success of virtual material design hinges on the question whether the chosen model provides a good fit to the considered dataset. If the model is too simple, then important characteristics of the dataset are missed. Hence, also the functional properties of the virtual realizations may differ substantially from the actual materials. On the other hand, extremely detailed models are at the risk of over-fitting the data. That is, realizations of such models may mimic certain particularities of the specific sample, and thereby fail to reflect accurately the true variability of microstructures seen in real data. In this setting, statistical hypothesis tests provide a scientifically sound methodology to decide whether at a given significance level a considered dataset is consistent with a null hypothesis formulated in advance.

%
%TEST STATISTIC
%
The design of versatile test statistics is the key challenge in devising powerful statistical hypothesis tests for data involving stochastic morphologies. The difficulty lies in finding appropriate numerical characteristics that encapsulate the relevant complex topology of the dataset.
In particular, different model classes from stochastic geometry, e.g., point processes, particle processes, line processes or random tessellations, come with different natural statistics. Here, we will restrict attention to random tessellations which are widely established as models for cellular and polycrystalline materials \cite{Bourne,Redenbach2009}.
The first natural step is to consider geometric properties of the tessellation cells such as edge lengths or cell surface areas  \cite{evolv}. However, for detecting subtle morphological differences, more refined characteristics may deliver additional value. For instance, Arns et al. \cite{arns} compute Minkowski functionals of morphologies resulting through continuous dilation of the original dataset.

%
%TDA
%
More recently, topological data analysis (TDA) has emerged as a promising new idea to detect subtle differences in the shapes of complex morphologies occurring in materials science \cite{tda_pnas,tda_nature}. Similarly as in \cite{arns}, this novel approach studies the evolution of the morphology resulting from suitable dilation operations. Instead of the Minkowski functionals, TDA uses the toolkit of persistent homology from algebraic topology to detect structure changes after dilation.  Hence, persistent homology is a characteristic, which by its design is a promising candidate when it comes to distinguishing structures with respect to their topological properties rather than their geometric ones.

%
%TESS-ADAPTED
%
To date, the widespread use of TDA in materials science is severely restricted since there is no tessellation-adapted version of the persistence diagram. More precisely, while the \v Cech and Vietoris-Rips filtrations are highly popular tools to compute the persistence diagram on point clouds, these filtrations do not take into account the rich adjacency structure inherent in spatial random tessellations. The key contribution of our work is to propose two specific ways to devise a variant of the persistence diagram that is capable of reflecting the information that the vertices belong to a random spatial tessellation. This will lead to the \emph{edge-based persistence diagram} and the \emph{$M$-localized persistence diagram} that we discuss in further detail in Section \ref{mod_sec} below.
To illustrate that TDA-based testing of spatial tessellation models is feasible, we carry out an analysis in three directions.

%
%RESULTS
%
First, we show that the persistent Betti numbers derived from the new persistence diagrams are asymptotically normal in large domains. We stress that we derive this asymptotic normality in the form of a functional central limit theorem. This functional formulation is an decisive advantage as it allows us to consider not only the individual persistent Betti numbers but any test statistics that results as a continuous transformation from the persistence diagram. This is a key insight since it guarantees that for large sampling windows, asymptotically exact hypothesis tests can be constructed from the knowledge of means and variances in the null model. In the context of spatial random tessellations, arguably the most fundamental null model is that of a \emph{Poisson-Voronoi tessellation}, corresponding to a tessellation obtained from a randomly scattered collection of indistinguishable cell centers. However, due to the complex physical phenomena governing the organization of real microstructures, the Poisson-Voronoi tessellation is rarely considered as a serious contender for an appropriate null model. Therefore, we establish the asymptotic normality in the framework of Laguerre tessellations with generators forming a Gibbsian point process. 
This enables modelling of substantial variations in the tessellation cells as well as interactions between the cell centers. 

Second, in a simulation study, we illustrate that the asymptotic normality becomes already accurate for moderately large sample sizes. We also give first indications of the power of the TDA-based test statistics and compare the testing power of these new test statistics to more elementary alternatives. Third, we apply the testing methodology to a specific dataset of an open cell foam from materials science.

We stress that the goal of our investigation is to develop a TDA-based framework for analyzing tessellations that 1) has a rigorous statistical foundation, 2) is validated in a simulation study, and 3) is applied to a challenging dataset from materials science. In particular, considered in itself, each of the three parts contains possibilities for further investigations. However, to strengthen the coherence between the different parts, we did not work out such refinements although they would be logical when considered in isolation within one of the parts. We believe that after having made the first steps in the present work, these could be exciting avenues for future research. We will elaborate on these possibilities in Section \ref{sec:cr}.

\section{Model and main results}
\label{mod_sec}

%MOT
We develop goodness of fit tests for random spatial tessellations observed in the sampling window $W_n := [0, n]^p$, $p \ge 2$. Due to their central place in materials science, we focus on the pivotal model of Laguerre tessellations, which include Voronoi tessellations as a special case \cite{Redenbach2009}.  Laguerre tessellations are obtained from a marked point process of cell centers through a deterministic construction rule.

%DEF
More precisely, we assume that $\XX_n := \XX \cap W_n$, where $\XX = \big\{\cX_i\big\}_{i \ge 1} = \{(X_i, R_i)\}_{i \ge 1}$ is a stationary marked point process, with locations $X_i \in \R^p$, and marked by radii $R_i > 0$. Then, we let
$$C\big(\cX_i, \XX_n\big) := \big\{y \in \R^p\co |y - X_i|^2 - R_i^2 \le \min_{(x, r) \in \XX_n}(|y - x|^2 - r^2)\big\}$$
denote the Laguerre cell associated with $\cX_i$. The collection $\Xi_n := \Xi(\XX_n) := \big\{C\big(\cX_i, \XX_n)\big\}_{X_i \in \XX_n}$ of all Laguerre cells defines the Laguerre tessellation, and we write $\Xi_n^{(q)}$ for the family of all $q$-faces of $\Xi_n$.

Another approach would be to build the tessellation on the entire stationary point process $\XX$, and then only restrict to the cells whose generator is contained in the window $W_n$. Although conceptually appealing, this approach has the disadvantage that the existing asymptotic theory for functionals on Gibbs point processes that will be used later is developed for functionals computed on $\XX_n$ \cite{gibbs}. It seems plausible that under suitable additional conditions on the process and the functional, the two approaches should become equivalent. However, already the verification of the standard conditions in \cite{gibbs} is delicate. For the sake of attaining a more accessible presentation, we therefore decided to refrain from establishing this equivalence in the present work.

\subsection{Edge-based $M$-bounded persistent Betti numbers}
\label{edge_exc}
%THICK
%FACE CLOSE
The permeability of an open cell foam is influenced by the size distribution of the windows between the foam cells \cite{Foehst,WESTHOFF20181}. Hence, when modeling the foam by a random tessellation, the distribution of face sizes has to be controlled. Suitable geometric characteristics are the area of faces or the \emph{face-inradius}. Here, we propose to proceed in the vein of \cite{arns} and consider the latter.

Fixing the dimension $p = 3$ and starting from the tessellation edges, the face inradius corresponds to the dilation radius at which a face gets completely covered by the edges forming its boundary. More precisely, we define the \emph{inradius} $\rfc(f)$ of a convex polygon $f$ bounded by the edges $\{e_1, \dots, e_k\}$ as
\begin{align}
	\label{ebf_eq}
	\rfc(f) := \max_{x \in f}\dist\big(x, \bigcup_{j \le k}e_{j}\big)
\end{align}
i.e., as the first level $r \ge 0$ where the $r$-dilated edge set covers the entire face. Due to edge effects, faces that are close to the boundary of the sampling window may be highly skewed and could therefore potentially influence the test statistics significantly.

To counter such effects, we will restrict our attention to faces with a bounded eccentricity: for a 2-face $f\in \Xi_n^{(2)}$ {shared by the cells} centered at $P_1$ and $P_2$, we define $\ms{ecc}(f) := \max_{y \in f, \, i \in\{1, 2\}} |y - P_i|$. This ensures that both the face is not too far from the centers, and also that its diameter is small.
Hence, for Voronoi tessellations, the eccentricity includes information on both, the size of a face and its distance to the generator points of the adjacent cells. For Laguerre tessellations based on independent radii  with a highly skewed distribution the interpretation is more involved. Indeed, for generators with small radius, it may happen that the generating point is not contained in the cell. Thus, in the Laguerre setting, the eccentricity also takes into account the degree of cell inhomogeneity.

Then, we introduce the \emph{edge-based persistent Betti number}
$$\b_n^{\ms e, M, s} := \#\big\{f \in \Xi_n^{(2)}\co \rfc(f) > s,\, \ms{ecc}(f) \le M\big\}$$
as the number of faces with inradius at least $s$ and eccentricity at most $M$.

\subsection{$M$-localized persistent Betti numbers}
\label{vert_exc}
\begin{figure}
	\begin{tikzpicture}[scale=0.7]
	\fill[black!10!white] (1.46, 0.10) circle (0.5cm);
\fill[black!10!white] (-0.2, 1.05) circle (0.5cm);
\fill[black!10!white] (-1.66, 0.00) circle (0.5cm);
\fill[black!10!white] (-0.00, -1.05) circle (0.5cm);
\draw[dotted] (1.46, 0.10) circle (0.5cm);
\draw[dotted] (-0.2, 1.05) circle (0.5cm);
\draw[dotted] (-1.66, 0.00) circle (0.5cm);
\draw[dotted] (-0.00, -1.05) circle (0.5cm);
\fill (1.46, 0.10) circle (2pt);
\fill (-0.2, 1.05) circle (2pt);
\fill (-1.66, 0.00) circle (2pt);
\fill (-0.00, -1.05) circle (2pt);
\end{tikzpicture}\begin{tikzpicture}[scale=0.7]
	\fill[black!10!white] (1.46, 0.10) circle (0.98cm);
\fill[black!10!white] (-0.2, 1.05) circle (0.98cm);
\fill[black!10!white] (-1.66, 0.00) circle (0.98cm);
\fill[black!10!white] (-0.00, -1.05) circle (0.98cm);
\draw[dotted] (1.46, 0.10) circle (0.98cm);
\draw[dotted] (-0.2, 1.05) circle (0.98cm);
\draw[dotted] (-1.66, 0.00) circle (0.98cm);
\draw[dotted] (-0.00, -1.05) circle (0.98cm);
\fill (1.46, 0.10) circle (2pt);
\fill (-0.2, 1.05) circle (2pt);
\fill (-1.66, 0.00) circle (2pt);
\fill (-0.00, -1.05) circle (2pt);
\draw[blue, ultra thick] (-1.66, 0.00)--(-0.2, 1.05);
\draw[blue, ultra thick] (-1.66, 0.00)--(-0.00, -1.05);
\draw[blue, ultra thick] (-0.2, 1.05)--(1.46, 0.10);
\draw[blue, ultra thick] (-0.00, -1.05)--(1.46, 0.10);
\end{tikzpicture}\begin{tikzpicture}[scale=0.7]
	\fill[black!10!white] (1.46, 0.10) circle (1.05cm);
\fill[black!10!white] (-0.2, 1.05) circle (1.05cm);
\fill[black!10!white] (-1.66, 0.00) circle (1.05cm);
\fill[black!10!white] (-0.00, -1.05) circle (1.05cm);
\draw[dotted] (1.46, 0.10) circle (1.05cm);
\draw[dotted] (-0.2, 1.05) circle (1.05cm);
\draw[dotted] (-1.66, 0.00) circle (1.05cm);
\draw[dotted] (-0.00, -1.05) circle (1.05cm);
\fill (1.46, 0.10) circle (2pt);
\fill (-0.2, 1.05) circle (2pt);
\fill (-1.66, 0.00) circle (2pt);
\fill (-0.00, -1.05) circle (2pt);
\draw[blue, ultra thick] (-1.66, 0.00)--(-0.2, 1.05);
\draw[blue, ultra thick] (-1.66, 0.00)--(-0.00, -1.05);
\draw[blue, ultra thick] (-0.2, 1.05)--(1.46, 0.10);
\draw[blue, ultra thick] (-0.00, -1.05)--(1.46, 0.10);
\draw[blue, ultra thick] (-0.00, -1.05)--(-0.2, 1.05);
\end{tikzpicture}\begin{tikzpicture}[scale=0.7]
	\fill[black!10!white] (1.46, 0.10) circle (1.16cm);
\fill[black!10!white] (-0.2, 1.05) circle (1.16cm);
\fill[black!10!white] (-1.66, 0.00) circle (1.16cm);
\fill[black!10!white] (-0.00, -1.05) circle (1.16cm);
	\fill[blue!30!white, opacity=0.8] (-1.66, 0.00)--(-0.2, 1.05)--(1.46, 0.10)--(-0.00, -1.05);
\draw[dotted] (1.46, 0.10) circle (1.16cm);
\draw[dotted] (-0.2, 1.05) circle (1.16cm);
\draw[dotted] (-1.66, 0.00) circle (1.16cm);
\draw[dotted] (-0.00, -1.05) circle (1.16cm);
\fill (1.46, 0.10) circle (2pt);
\fill (-0.2, 1.05) circle (2pt);
\fill (-1.66, 0.00) circle (2pt);
\fill (-0.00, -1.05) circle (2pt);

\draw[blue, ultra thick] (-1.66, 0.00)--(-0.2, 1.05);
\draw[blue, ultra thick] (-1.66, 0.00)--(-0.00, -1.05);
\draw[blue, ultra thick] (-0.2, 1.05)--(1.46, 0.10);
\draw[blue, ultra thick] (-0.00, -1.05)--(1.46, 0.10);
\draw[blue, ultra thick] (-0.00, -1.05)--(-0.2, 1.05);
\end{tikzpicture}\\
	\begin{tikzpicture}[scale=0.7]
	\fill[black!10!white] (1.46, 0.10) circle (0.5cm);
\fill[black!10!white] (-0.2, 1.05) circle (0.5cm);
\fill[black!10!white] (-1.66, 0.00) circle (0.5cm);
\fill[black!10!white] (-0.00, -1.05) circle (0.5cm);
\draw[dotted] (1.46, 0.10) circle (0.5cm);
\draw[dotted] (-0.2, 1.05) circle (0.5cm);
\draw[dotted] (-1.66, 0.00) circle (0.5cm);
\draw[dotted] (-0.00, -1.05) circle (0.5cm);
\fill (1.46, 0.10) circle (2pt);
\fill (-0.2, 1.05) circle (2pt);
\fill (-1.66, 0.00) circle (2pt);
\fill (-0.00, -1.05) circle (2pt);
\draw[dashed] (-1.66, 0.00)--(-0.2, 1.05);
\draw[dashed] (-1.66, 0.00)--(-0.00, -1.05);
\draw[dashed] (-0.2, 1.05)--(1.46, 0.10);
\draw[dashed] (-0.00, -1.05)--(1.46, 0.10);
\end{tikzpicture}\begin{tikzpicture}[scale=0.7]
	\fill[black!10!white] (1.46, 0.10) circle (0.98cm);
\fill[black!10!white] (-0.2, 1.05) circle (0.98cm);
\fill[black!10!white] (-1.66, 0.00) circle (0.98cm);
\fill[black!10!white] (-0.00, -1.05) circle (0.98cm);
\draw[dotted] (1.46, 0.10) circle (0.98cm);
\draw[dotted] (-0.2, 1.05) circle (0.98cm);
\draw[dotted] (-1.66, 0.00) circle (0.98cm);
\draw[dotted] (-0.00, -1.05) circle (0.98cm);
\fill (1.46, 0.10) circle (2pt);
\fill (-0.2, 1.05) circle (2pt);
\fill (-1.66, 0.00) circle (2pt);
\fill (-0.00, -1.05) circle (2pt);
\draw[blue, ultra thick] (-1.66, 0.00)--(-0.2, 1.05);
\draw[blue, ultra thick] (-1.66, 0.00)--(-0.00, -1.05);
\draw[blue, ultra thick] (-0.2, 1.05)--(1.46, 0.10);
\draw[blue, ultra thick] (-0.00, -1.05)--(1.46, 0.10);
\end{tikzpicture}\begin{tikzpicture}[scale=0.7]
	\fill[black!10!white] (1.46, 0.10) circle (1.05cm);
\fill[black!10!white] (-0.2, 1.05) circle (1.05cm);
\fill[black!10!white] (-1.66, 0.00) circle (1.05cm);
\fill[black!10!white] (-0.00, -1.05) circle (1.05cm);
\draw[dotted] (1.46, 0.10) circle (1.05cm);
\draw[dotted] (-0.2, 1.05) circle (1.05cm);
\draw[dotted] (-1.66, 0.00) circle (1.05cm);
\draw[dotted] (-0.00, -1.05) circle (1.05cm);
\fill (1.46, 0.10) circle (2pt);
\fill (-0.2, 1.05) circle (2pt);
\fill (-1.66, 0.00) circle (2pt);
\fill (-0.00, -1.05) circle (2pt);
\draw[blue, ultra thick] (-1.66, 0.00)--(-0.2, 1.05);
\draw[blue, ultra thick] (-1.66, 0.00)--(-0.00, -1.05);
\draw[blue, ultra thick] (-0.2, 1.05)--(1.46, 0.10);
\draw[blue, ultra thick] (-0.00, -1.05)--(1.46, 0.10);
\end{tikzpicture}\begin{tikzpicture}[scale=0.7]
	\fill[black!10!white] (1.46, 0.10) circle (1.16cm);
\fill[black!10!white] (-0.2, 1.05) circle (1.16cm);
\fill[black!10!white] (-1.66, 0.00) circle (1.16cm);
\fill[black!10!white] (-0.00, -1.05) circle (1.16cm);
	\fill[blue!30!white, opacity=0.8] (-1.66, 0.00)--(-0.2, 1.05)--(1.46, 0.10)--(-0.00, -1.05);
\draw[dotted] (1.46, 0.10) circle (1.16cm);
\draw[dotted] (-0.2, 1.05) circle (1.16cm);
\draw[dotted] (-1.66, 0.00) circle (1.16cm);
\draw[dotted] (-0.00, -1.05) circle (1.16cm);
\fill (1.46, 0.10) circle (2pt);
\fill (-0.2, 1.05) circle (2pt);
\fill (-1.66, 0.00) circle (2pt);
\fill (-0.00, -1.05) circle (2pt);

\draw[blue, ultra thick] (-1.66, 0.00)--(-0.2, 1.05);
\draw[blue, ultra thick] (-1.66, 0.00)--(-0.00, -1.05);
\draw[blue, ultra thick] (-0.2, 1.05)--(1.46, 0.10);
\draw[blue, ultra thick] (-0.00, -1.05)--(1.46, 0.10);
\end{tikzpicture}
	\caption{Top row: Snapshots at four time points illustrating the evolution of the standard \v Cech filtration for a  cloud of four points. At time point 1, all four  points are in separate connected components. At time point 2, a loop has appeared, thereby leading to the birth of a 1-feature. At time point 3 a second 1-feature is born. At time point 4 both 1-features die. Bottom row: Corresponding evolution of the tessellation-adapted filtration, where the four points are the vertices of a cell {(in dashed lines)}. At time points 1 and 2 there is no difference to the standard \v Cech filtration. However, at time point 3 no new feature is born since the connecting edge is not part of the tessellation. At time point 4, the feature born at time point 1 dies.}
	\label{filt_fig}
\end{figure}
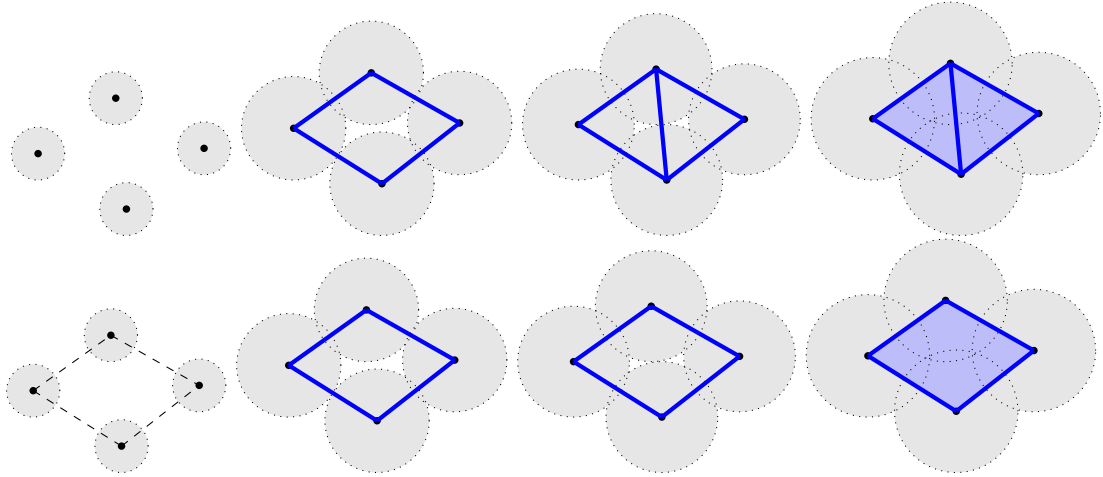

Face-inradii are very natural univariate functional characteristics associated with a tessellation. However, TDA draws its appeal from the ability to extract far more subtle bivariate quantities describing the times when topological features appear, and when they disappear again.
This machinery relies on a suitable notion of filtrations of topological spaces. In the analysis of point patterns, the \v Cech- and the Vietoris-Rips filtrations have emerged as universally accepted choices. However, when dealing with tessellations there is no clear suggestion from the literature that one could build upon. In order to respect the inherent structure of the tessellation, we propose tessellation-adapted filtrations where entire lower-dimensional faces are added at filtration times that are given by their circumradii. In other words, a $q$-face $f \in \Xi_n^{(q)}$ defined by vertices $P_0, \dots, P_m$ belongs to the filtration at a level $s > 0$ if and only if $\rci(f) \le s$, where the  \emph{circumradius} $\rci(f)$ of $f$ is defined by
$$\rci(f) := \min_{y \in \R^p}\max_{i \le m} |y - P_i|.$$
 The circumradius of a $q$-face is always determined by the circumradii of its $q$-simplices in the sense that there exist $P_{j_0}, \dots, P_{j_q}$ with $\rci(f) = \rci(\{P_{j_0}, \dots, P_{j_q}\})$. We will refer to this as the \emph{tessellation-adapted filtration} associated with the faces of the tessellation. Figure \ref{filt_fig} illustrates the difference between the tessellation-adapted and the standard \v Cech filtration.

%
%
%

%
%PERS DIAG
%
Next, we construct the \emph{persistence diagram} $\{(B_{i, n, q}, D_{i, n, q})\}_{i \ge 1}$ of birth times $B_{i, n, q}$ and death times $D_{i, n, q}$ of the $q$-features. For $p=3$, a $0$-feature is a connected component, a $1$-feature is a circular hole (tunnel), and a $2$-feature is a void (connected component of the complement). Based on the persistence diagram, we define the \emph{persistent Betti number} for $b, d > 0$ as
\begin{align}
	\label{bet_eq}
	\b_{n, q}^{b, d} := \#\{i \ge 1\co B_{i, n, q} \le b, D_{i, n, q} \ge d\},
\end{align}
the number of $q$-features born before time $b > 0$ and living past time $d > 0$.
\medskip

%
%SCALAR CLT
%
Since it is highly delicate to derive the persistence diagram from a given point cloud, we need to implement some approximation and truncation steps in order to establish the asymptotic normality rigorously.  First, we want to avoid very elongated faces with a small inradius, and therefore write $\Xi_{n, M}^{',(q)}$ for the family of all $q$-faces in $\Xi_n$ of eccentricity at most $M$ and whose inradius is at least $M^{-1}$ for some (large) $M > 0$. An additional major difficulty arises through long-range correlations induced by very large features spanning over large parts of the sampling window. To deal with this issue, we restrict our attention to \emph{$M$-localized Betti numbers}, which we now discuss in further detail. Loosely speaking, $M$-localized persistent Betti numbers first determine locally the contribution of features involving a $q$-face and then aggregate these contributions over all faces in the entire window. More precisely, fixing $M>0$, $b \le d$, for each $f \in \Xi_{n, M}^{', (q)}$, we consider all cells contained in the sub-window $z(f) + [-M, M]^p$ centered at the centroid $z(f)$ of the face $f$.

When computing the persistent Betti numbers, we rely on an incremental algorithm. We now explain loosely how this algorithm can be used to compute the birth times and death times of features, referring the reader to \cite[Algorithm 11]{yvinec} for a detailed exposition. We start from a completely empty space, and then reconstruct the tessellation by adding the faces one after another in the order of increasing circumradius. When adding a face, one of two cases may occur. Either the addition of the face kills one of the existing features, in which case it is called \emph{negative}. For instance, in dimension $p = 3$, the addition of a face can cause a loop to be contractible. Otherwise, the face is called \emph{positive} as it will give rise to a new feature. In $p = 3$, this corresponds to the creation of a 3D cavity. If the face $f$ is negative, then we write $b(f) $ for the birth time and $d(f)= \rci(f)$ for the death time of the feature killed by the addition of face $f$.   In the case where the feature-death is caused by the simultaneous addition of multiple faces, we can formally attach the feature to the face with the left-most circumcenter. Note that this case occurs with probability 0 for the considered random tessellations. The \emph{$M$-localized persistent Betti number $\b_{n, q}^{M, b, d}$} is then defined as the total number of such features
	$$\b_{n, q}^{M, b, d} := \#\big\{f \in \Xi_{n, M}^{',(q)}\co b(f) \le b,\, d(f) \ge d\big\}.$$

\medskip

\subsection{Definition of the generating point process and conditions for the main result}
\label{cond_sec}

Since tessellations appearing in materials science typically feature highly complex interactions between different cells, it is essential to work with classes of marked point processes that are flexible enough to reproduce such interaction patterns. Since Gibbsian point processes have been successfully used in a variety of application scenarios with complex dependencies, our investigations will be focused on this class of point processes.

%
%PRIORI
%
First, we recall the definition of Gibbs point processes from \cite[Equation (1.5)]{gibbs}. Let $\PP_{\t, \Q}$ denote an independently marked Poisson point process on $\R^p$ with intensity $\t > 0$ and mark  distribution $\Q$  whose support will be assumed to be contained in a bounded interval of the form $[0, \mu]$ for some $\mu > 0$.

%
%GIBBS DEF
%
Second, let $\Psi$ be a translation- and rotation-invariant nonnegative functional defined on finite point configurations, which we henceforth refer to as a \emph{Hamiltonian}.
For $\a > 0$, a stationary marked point process $\XX^{\a\Psi}$ is called \emph{Gibbs point process with Hamiltonian $\a\Psi$}, if for every bounded open $D \su \Rm :=\R^p \ti [0, \mu]$, the Radon-Nikodym density of the marked point process $\XX_D = \XX \cap D$ with respect to the Poisson point process $\PP_{\t, \Q} \cap D$ is proportional to $e^{-\a \Psi(\cdot )}$.
That is, writing $\mc L$ for the law of a point process, we have
$$\f{\d \mc L(\XX_D)}{\d \mc L(\PP_{\t, \Q} \cap D)}(\vp) = \f1{Z_{\a, D}}{e^{-\a \Psi(\vp )}},$$
where $\vp \su D$ is a finite configuration and $Z_{\a, D} := \E[\exp(-\a \Psi(\PP_{\t, \Q} \cap D))]$.
In the following, we set $\XX = \XX^{\a\Psi}$. If $\a$ and $\Psi$ are fixed, we will also use the abbreviation $\XX_n := \XX_{W_n \ti [0, \mu]}^{\a \Psi}$.

In the following, we will prove asymptotic normality of test statistics based on the edge-based and $M$-localized persistent Betti numbers. Recently, establishing functional CLTs for Gibbs point patterns has become a vigorous research topic in spatial statistics. For instance, \cite{svane2} establishes asymptotic normality for Ripley's $K$-function. This quantity is an important characteristic in spatial statistics but only takes into account the distances between point pairs of the Gibbs process. The setting in the present paper is substantially more involved since we deal with topological functionals of 3D tessellations where the Gibbs point patterns are not the vertices but the generators of the cells. On the conceptual side, one of the major contributions of our work is to explain how the delicate moment bounds can be established in this complex setting. Similarly as in \cite{svane2}, our proof crucially relies on the techniques from \cite{gibbs,gibbsCLT}. Although, formally the framework there is stated for point processes without marks, the remark following \cite[Theorem 2.4]{gibbs} stresses that the proof of the asymptotic normality also extends seamlessly to marked Gibbs point processes.

%
%HAMILTONIAN
%
The asymptotic normality in \cite{gibbs,gibbsCLT} is established for a specific class $\bs \Psi^*$ of translation- and rotation-invariant nonnegative  Hamiltonians $\Psi$ defined on finite point configurations $\vp \su \Rm$: a Hamiltonian $\Psi$ belongs to the class $\bs \Psi^*$ if   i) $\Psi(\vp) \le \Psi(\vp')$ for $\vp \su \vp'$, ii) $\Psi(\vp) < \ff$ whenever $\# \vp \le 1$, and iii) $\Psi$ has a finite range $r_0 > 0$. We refer the reader to \cite[Section 1.1.1]{gibbsCLT} for a more detailed discussion. For the rest of the present paper, the main implication is that the Hamiltonian of the Strauss process considered in Section \ref{sim_sec} below belongs to the class $\bs \Psi^*$.

\smallskip

%
%TES CON
%

 \subsection{Asymptotic normality}
	 \label{res_sec}
 After having introduced the conditions in Section \ref{cond_sec}, we now state asymptotic normality of test statistics computed from edge-based and $M$-localized persistent Betti numbers. First, we consider statistics, which concern persistent Betti numbers for fixed values of the levels.
 To state our results precisely, we require an upper bound for the Poisson intensity which is given by
 $$\t_0(\Psi) :=  |B_1(o)|^{-1}\exp(\a m_0^\Psi)(r_0 + 1)^{-p},$$
 where
$m_0 ^\Psi := \inf_{\vp \text{ loc.~finite}} \big(\Psi(\vp \cup \{o \}) - \Psi(\vp)\big).$

%
%SCLT
%
\bet[Asymptotic normality at a fixed level]
\label{clt}
Let $\XX = \XX^{\a\Psi}$ be a marked Gibbs point process with Hamiltonian $\a\Psi$ of type $\bs \Psi^*$, and $\t < \t_0(\Psi)$.
\been
\im
Let $s > 0$. Then,
$$\f{\b_{ n}^{\ms e, M, s} - \E[\b_{ n}^{\ms e, M, s}]}{\sqrt{|W_n|}}$$
converges in distribution to a mean-zero normal random variable.
\im
Let  $d \ge b \ge 0$. Then,
$$\f{\b_{n, p - 2}^{M, b, d} - \E[\b_{n, p - 2}^{M, b, d}]}{\sqrt{|W_n|}}$$
converges in distribution to a mean-zero normal random variable.
\enen
\ent
It is possible to provide an integral expression for the variance of the limiting random variable, see \cite[Identity 2.11]{gibbs}. However, since this integral is not amenable to an evaluation in closed form, we refrain from reproducing the precise expression. In the simulation study in Section \ref{sim_sec} below, we determine the variance under the null model through simulations.

Finally, we extend the fixed-level asymptotic normality from Theorem \ref{clt} to asymptotic normality on a functional level. That is, we consider $\b_{ n}^{\ms e, M, s}$ as a function in $s$ and $\b_{n, p - 1}^{M, b, d}$ as a function in $b, d$. While for the fixed-level asymptotic normality, the main challenge is to control the long-range spatial dependence, the extension to the functional level also relies on continuity properties when varying the levels at which the persistence Betti numbers are computed.  Establishing the required continuity properties is already delicate when dealing with point clouds described by a  Poisson point process, and a direct extension to Gibbs-Laguerre tessellations seems out of reach. We now explain that these problems can be avoided by incorporating additional sources of noise into the original models.

More precisely, motivated by applications such as foams we assume that the tessellation edges are not just purely 1-dimensional links but come with a certain random thickness. That is, we assume that the process of edges $\{e_i\}_{i \ge 1}$ is independently marked with positive thicknesses $\{\r_i\}_{i \ge 1}$. This model extension helps to ensure a higher degree of distributional smoothness on the conceptual side and also moves the model closer to real structures found in materials science.

In the case of the $M$-localized persistent Betti numbers, we assume that  extracting the location of the tessellation vertices from data comes with a certain measurement error \cite{cheng}. We incorporate this constraint by considering iid noise vectors $\{\eta_i\}_{i \ge 1}$ that are uniformly distributed in a ball $B_{h_0}$ for some small $h_0 > 0$. Then, we set the measured locations as $P_i + \eta_i$, which will replace $P_i$ when determining the circumradius in the tessellation-adapted filtration.

%
%FCLT
%
\bet[Functional asymptotic normality]
\label{fclt}
Let $S > 0$ and $\XX = \XX^{\a\Psi}$ be a marked Gibbs point process with Hamiltonian $\a\Psi$ of type $\bs \Psi^*$ and $\t < \t_0(\Psi)$.
\been
\im When considered as a process in $s \le S$, the process 
$$\big(|W_n|^{-1/2}(\b_n^{\ms e, s} - \E[\b_n^{\ms e, s}])\big)_{s \le S}$$
converges in the Skorokhod topology to a centered Gaussian process.
\im When considered as a process in $b, d \le S$, the process
$$\big(|W_n|^{-1/2}(\b_{n, p - 2}^{M, b, d} - \E[\b_{n,  p - 2}^{M, b, d}])\big)_{s \le S}$$
converges in the Skorokhod topology to a centered Gaussian process.
\enen
\ent

We note that for the $M$-localized persistent Betti numbers, Theorems \ref{clt} and \ref{fclt} could also be shown for $\b_{n, p - 1}$ instead of $\b_{n, p - 2}$ when associating a feature with the face  causing its birth (instead of the death).  Moreover, the restriction to $M$-localized features can be avoided when replacing the persistent Betti numbers by the Euler characteristic curve. Although the latter is a coarser characteristic, it may still deliver interesting topological insights.

\section{Simulation study}
\label{sim_sec}
In this section, we present a simulation study analyzing goodness-of-fit tests derived from the asymptotic normality established in Theorems \ref{clt} and \ref{fclt}. To that end, we consider six tessellation models on the 3D unit cube $[0, 1]^3$, where to compensate for edge effects, we work with periodic boundary conditions. The first three models correspond to Voronoi tessellations formed on configurations of 300 generators in the unit cube. In the simplest model, the 300 generators are scattered independently at random in $[0, 1]^3$. In the second model, the 300 generator points are distributed according to a Strauss process, which is a Gibbs point process featuring a repulsive interaction between points. While the general limit theory developed in Section \ref{mod_sec} concerns Gibbs point processes with a variable number of points, for the simulation study, we have decided to rely on a variant where the number of points is held constant at 300. This model choice enables a cleaner comparison of the refined topological properties of the different tessellation models. We recall that the Strauss process is a Gibbsian point process whose unnormalized density with respect to a Binomial point process is given by $f(\xx) := \gamma^{s_{r_0}(\xx)}$, for certain parameters $r_0 > 0$ and $\g \in [0, 1]$. Here,  $s_{r_0}(\xx)$ is the number of point pairs at distance at most $r_0$. In particular, the inverse temperature $\a$ from Section \ref{cond_sec} and the parameter $\g$ are related through $\a = -\log(\g)$. 

Although Strauss processes are a popular choice for modeling repulsive point patterns, they are most useful for scenarios with a relatively low regularity. Therefore, we now explain more precisely how to create more regular point patterns. One of the most popular approaches is to construct random sphere packings via collective rearrangement, e.g. the force-biased algorithm. Such algorithms first start from a configuration of potentially overlapping spheres. Then, colliding spheres are iteratively pushed into the vacant space, and this procedure is repeated until a hard-core configuration is obtained \cite{forceBiased}. Such point patterns form the third model considered in our simulation study. We now summarize again the three models and provide the specific parameter choices.

\been
\im A {\bf Binomial-Voronoi tessellation}; short $\ms{Bin}$-$\ms{Vor}$. The Voronoi tessellation is constructed on 300 generator points that are scattered independently uniformly in $[0, 1]^3$. 
\im A {\bf Voronoi tessellation} based on a {\bf Strauss process}; short $\ms{St}$-$\ms{Vor}$. The interaction parameter $\g$ of the Strauss process is set to $\g = 0.01$, and the interaction range $r_0$ is fixed as $r_0 = 0.14$. When interpreting the point configuration as a set of non-overlapping balls of radius $r_0/2=0.07$, this corresponds to a volume fraction of 43\%. We note that larger choices of $r_0$ caused serious convergence issues of the MCMC sampler.
\im A {\bf Voronoi tessellation} based on a {\bf force-biased sphere packing} with a constant sphere radius of $r=0.07816$ resulting in a volume fraction of 60\%; short $\ms{Fb}$-$\ms{Vor}$ .
\enen
As a further refinement of the above setting, we replace the simplistic Voronoi model by a more refined Laguerre model with random radii. We note that to achieve on average 300 cells in the Laguerre tessellations, a substantially higher initial number of generator points is needed in order to compensate for empty cells. More precisely, the initial number of the generators was chosen as 325 for the Binomial process and 324 for the Strauss process. Since the spheres in the force-biased algorithm are disjoint by construction, no correction is needed in this example.
\been
\im[(4)] A {\bf Laguerre tessellation} based on a {\bf Binomial point process}; short $\ms{Bin}$-$\ms{Lag}$. The volume of the balls associated with the Laguerre tessellation follows a lognormal-distribution with parameters $\mu= -6.3262$ and $\sigma= 0.47238$. This corresponds to a coefficient of variation of 0.5.
\im[(5)] A {\bf Laguerre tessellation} based on a {\bf Strauss process}; short $\ms{St}$-$\ms{Lag}$. The interaction parameter $\g$ of the Strauss process is set to $\g = 0.01$, and the interaction radius $r_0$ is fixed as $r_0 = 0.14$. The Laguerre radius distribution is the same as before.
\im[(6)] A {\bf Laguerre tessellation} based on a {\bf force-biased sphere packing}; short $\ms{Fb}$-$\ms{Vor}$. The Laguerre radius distribution is the same as before and the volume fraction is 60\%.

\begin{figure}[!htpb]
    \centering
	\includegraphics[width = .79\textwidth]{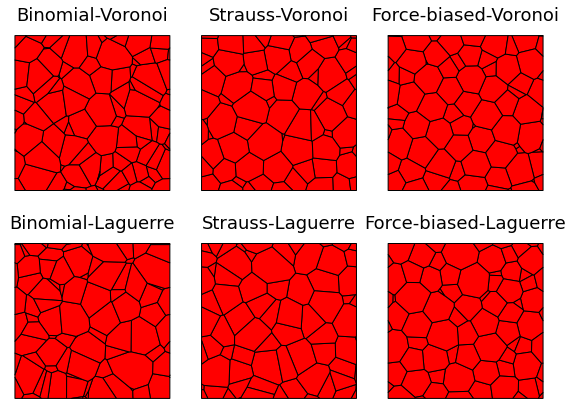}
	\caption{2D sections for realizations from each of the six tessellation models.}
    \label{sec_fig}
\end{figure}

In Figure \ref{sec_fig}, we illustrate 2D sections for realizations from each of the six tessellation models. 

\enen

\subsection{Section overview}

First, Section \ref{expl_sec} reveals how geometric and topological properties of the different models are reflected by suitable test statistics. Then, Section \ref{norm_sec} illustrates to which extent the approximate asymptotic normality is accurate in bounded sampling windows and for models that go beyond those covered by the sufficient conditions in Theorems \ref{clt} and \ref{fclt}. Finally, Section \ref{pow_sec} compares the power of goodness-of-fit tests for the statistics considered in Sections \ref{expl_sec} and \ref{norm_sec}.
\vspace{1cm}
%
%EXPL
%
\subsection{Exploratory analysis}
\label{expl_sec}

\begin{figure}[!ht]
    \centering
	\includegraphics[width = .99\textwidth]{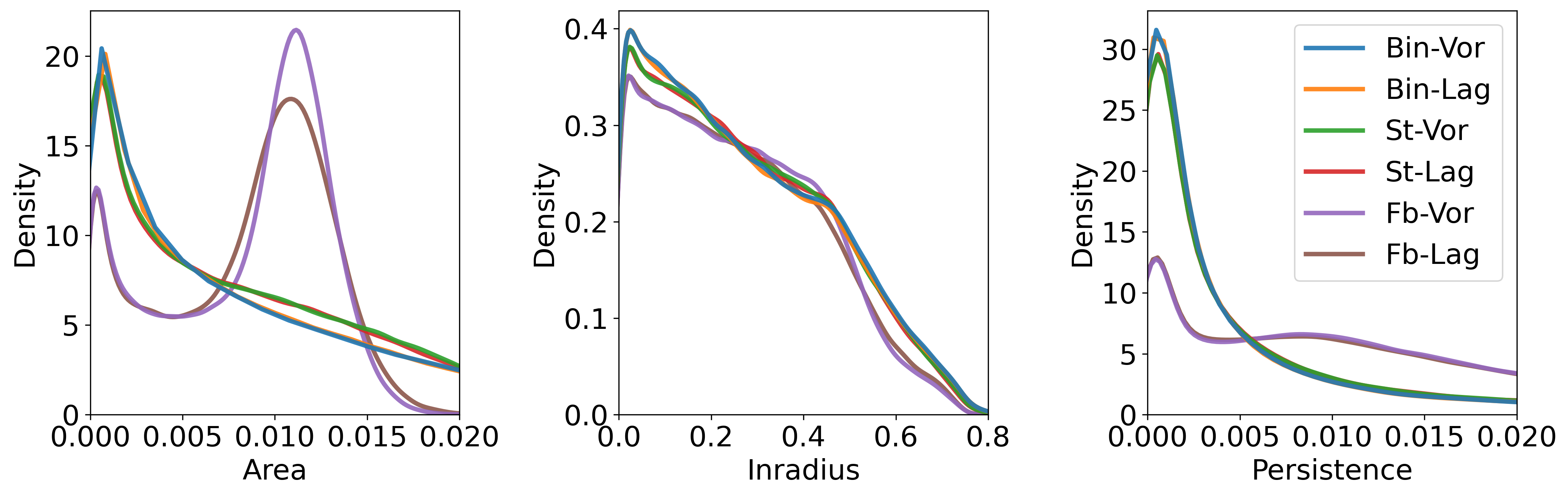}
	\caption{Densities of face areas (left), face inradii (middle) and persistences (right). }
    \label{kde_fig}
\end{figure}

In this section, we illustrate how different test statistics behave for the six tessellation models described in Section \ref{sim_sec}. We illustrate in the kernel-density estimates in Figure \ref{kde_fig}  three prototypical examples of model characteristics, namely \emph{face areas}, \emph{face inradii} and \emph{vertex-based persistences}.

To begin with, we consider the distribution of the face areas in the six models. Due to the higher variability, faces with a very small area occur far more often in the Binomial than in the Strauss models. This difference becomes even more pronounced when focussing on exceptionally large cells. Moreover, in the force-biased models the histograms have a clear peak corresponding to an area that is substantially bounded away from 0. Qualitatively, we can make similar observations for the persistence (i.e., lifetime) of features in the vertex-based filtration from Section \ref{vert_exc}. In particular, the kernel-density estimates in Figure \ref{kde_fig} illustrate that long lifetimes of features are common in the force-biased model. For the face-inradii from Section \ref{edge_exc}, the differences between the models are far less pronounced. Moreover, we observe that in all cases, the histograms corresponding to the Voronoi and to the Laguerre model are very similar.

%
%AS NORM
%
\subsection{Asymptotic normality}
\label{norm_sec}
The exploratory analysis in Section \ref{expl_sec} already revealed that suitable geometric and topological statistics can uncover substantial differences between considered models. In the next step, we use these summary statistics to devise goodness-of-fit tests. We have experimented with different forms of statistics and in general found the ones concentrating on the upper tails to have the highest testing power. Hence, we will use the following three statistics.

%STATS
\been
\im {\bf Face areas.} The statistic
$$\TA :=  \#\{f \in \Xi_n^{(2)}\co |f| > a_{\ms{Area}}\}$$
counts the number of tessellation faces $f$ with area $|f|$ exceeding a given threshold $a_{\ms{Area}} > 0$.
\im {\bf Face inradii.} The statistic
$$\TF :=  \#\{f \in \Xi_n^{(2)}\co \rfc(f) > a_{\ms{I}}\}$$
counts the number of tessellation faces $f$ with face-inradii time $\rfc(f)$ exceeding a given threshold $a_{\ms{I}} > 0$.
\im {\bf Persistence in tessellation-adapted filtrations.} The statistic
$$\TP :=  \#\{j \co D_j - B_j > a_{\ms{Pers}}\}$$
counts the number of features in the tessellation-adapted persistence diagram from Section \ref{vert_exc} with lifetime exceeding a given threshold $a_{\ms{Pers}} > 0$.
\enen

For selecting critical values of the tests, we rely on the approximate normality in large domains. We will first discuss why asymptotic normality is plausible from a conceptual view point (especially in light of Theorems \ref{clt} and \ref{fclt}), before providing numerical evidence that the normal approximation is already accurate in the considered sampling windows. 
 
\begin{figure}[!ht]
    \centering
	\includegraphics[width = .99\textwidth]{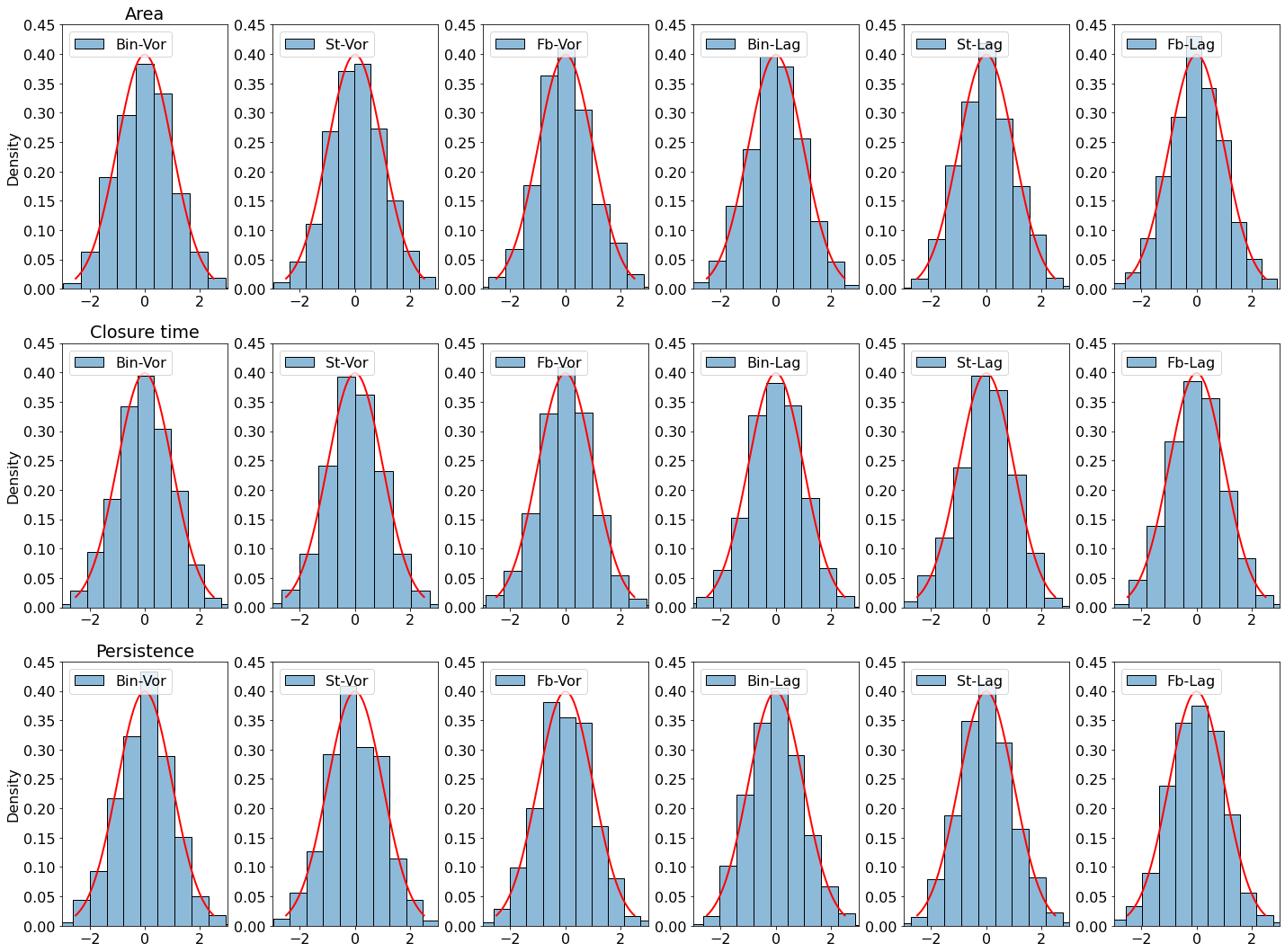}
	\caption{Standardized distributions of $\TA, \TF$, and $\TP$ (top, middle, and bottom). From left to right, each row shows the Binomial-Voronoi, Strauss-Voronoi, force-biased Voronoi, Binomial-Laguerre,  Strauss-Laguerre and force-biased Laguerre model.}
    \label{norm_fig}
\end{figure}

 First, we note that the number of tessellation faces whose inradius exceeds a fixed threshold corresponds to the edge-based persistent Betti numbers $\beta_n^{e,M,s}$, where to simplify the implementation, we have not enforced the technical restriction to faces with eccentricity at most $M$. As we will discuss in the following paragraph, although some of the technical assumptions needed in the proof of Theorem \ref{clt} are not satisfied in our simulation study, the numerical evidence provided in Figure \ref{norm_fig} illustrates convincingly that the asymptotic normality remains plausible in this more general setting. 
 
 For Poisson-Voronoi and Laguerre tessellations, \cite[Theorem 2.3]{flimmel} shows asymptotic normality of geometric functionals of cells subject to certain stabilization and moment conditions. In particular, the distribution function of the cell volume and the aggregated cell surface area are provided as specific examples. It is plausible that the methods can be extended to yield asymptotic normality of more complicated statistics on the face areas, and also when suitable Gibbs processes are used as generators. Again, our numerical results will make this belief credible. 
 
 Finally, we move to the $M$-localized persistent Betti numbers. Here, we recall that for fixed $d \ge b \ge 0$ the test statistic $\b_{n, 1}^{M, b, d}$ counts the number of 1-features that are born before time $b$ and die after time $d$. However, to recover the number of features exceeding a given lifetime a single test statistic $\b_{n, 1}^{M, b, d}$ does not suffice. This is similar to the situation considered in \cite{svane}. There, it was shown that the sum of all lifetimes is a continuous functional in $\b_{n, 1}^{M, b, d}$, when considered as a 2-parameter process in the entries $b, d\ge 0$. This stresses the importance of functional limit theorems (such as Theorem \ref{fclt}), since then the classical continuous mapping theorem makes it possible to derive the asymptotic normality for persistence-based test statistics that go beyond simple linear combinations of individual persistent Betti numbers.

In Theorems \ref{clt} and \ref{fclt}, the normality of the considered test statistics is established only in the limit of unboundedly large windows, this leaves the question to what extent the normal distribution is already appropriate for bounded sampling windows. Moreover, as mentioned above, some of the models described in Section \ref{sim_sec} do not satisfy all of the technical conditions that we imposed in order to be able to prove Theorems \ref{clt} and \ref{fclt} rigorously, namely:
\been
\im the number of cells is fixed in the Voronoi case.
\im the tessellation edges have no thickness;
\im the tessellation vertices are not perturbed;
\im the inequality $\t \le |B_1(o)|^{-1}\exp(\b m_0^\Psi)(r_0 + 1)^{-p}$ is violated;
\im the support of the mark distribution $\Q$ may be unbounded;
\im the restriction of considering only $M$-localized features is dropped;
\im the restrictions on $\ms{size}(f)$ are dropped;
\im the tessellation-adapted filtration is formed with respect to the alpha complex rather than the \v Cech complex.
\enen
We stress that we are confident that the asymptotic normality of Betti numbers holds under less restrictive conditions than the ones described in Theorems \ref{clt} and \ref{fclt}. As we will elaborate in further detail below, numerical evidence for this hypothesis is one of the major contributions from our simulation study.

%
%HIST FIG
%
In accordance with the tests to be derived in Section \ref{pow_sec}, we fix the thresholds $a_{\ms{Area}}$, $a_{\ms{I}}$, and $a_{\ms{Pers}}$ as the 40\%, 60\%, and 70\% quantiles of the corresponding test statistics under the considered null model. To fix the specific threshold levels, we used several pilot runs to determine which values have the potential to lead to promising goodness-of-fit tests.

Figure \ref{norm_fig} compares  the density of the standard normal distribution with the standardized test statistics based on 1,000 simulation runs. On a very general level, we draw the conclusion that the asymptotic normality is already reasonably accurate in moderately large windows.

%
%POWER
%
\subsection{Power analysis}
\label{pow_sec}

Now, we analyze the power of goodness-of-fit tests derived from the statistics $\TA$, $\TF$ and $\TP$ for the six models introduced in the beginning of Section \ref{sim_sec}. To that end, we proceed as follows. First, we generate 1,000 realizations from each of the six models to determine the distribution of the test statistics under the respective null hypotheses.

 \begin{table}[!htpb]
         \begin{center}
		 \caption{Rejection rates in percentage points for the goodness-of-fit tests based on  $\TA$, $\TF$ and $\TP$. Each row describes the null model that is used to compute the mean and variance for the test based on asymptotic normality; each column corresponds to the alternative for which the test statistic is computed. }
                 \label{pow_mark_tab}
		 \begin{tabular}{l|p{29pt}p{29pt}p{29pt}p{29pt}p{29pt}p{29pt}}
			 \multicolumn{7}{c}{$\TA$}\\
			 $\ms H_0\backslash \ms H_1$\hspace{-.15cm}            & $\ms {Bin}\text{-}\ms {Vor}$& $\ms {St}\text{-}\ms {Vor}$& $\ms {Fb}\text{-}\ms {Vor}$& $\ms{Bin}\text{-}\ms {Lag}$ & $\ms {St}\text{-}\ms {Lag}$ &$\ms {Fb}\text{-}\ms {Lag}$ \\
                         \hline

$\ms{Bin}\text{-}\ms{Vor}$&4.7&98.7&100.0&22.8&92.1&100.0\\
$\ms{St}\text{-}\ms{Vor}$&99.0&5.5&100.0&97.5&15.1&100.0\\
$\ms{Fb}\text{-}\ms{Vor}$&100.0&100.0&5.5&100.0&100.0&99.9\\
$\ms{Bin}\text{-}\ms{Lag}$&0.5&86.4&100.0&6.4&74.3&100.0\\
$\ms{St}\text{-}\ms{Lag}$&89.7&1.3&100.0&86.9&3.4&100.0\\
$\ms{Fb}\text{-}\ms{Lag}$&100.0&100.0&97.2&100.0&100.0&5.5

		 \end{tabular}\,

                 \begin{tabular}{l|p{29pt}p{29pt}p{29pt}p{29pt}p{29pt}p{29pt}}
			 \multicolumn{7}{c}{$\TF$}\\
			 $\ms H_0\backslash \ms H_1$\hspace{-.15cm}            & $\ms {Bin}\text{-}\ms{Vor}$& $\ms{St}\text{-}\ms{Vor}$& $\ms{Fb}\text{-}\ms{Vor}$& $\ms {Bin}\text{-}\ms{Lag}$ & $\ms{St}\text{-}\ms{Lag}$ &$\ms{Fb}\text{-}\ms{Lag}$ \\
                         \hline

$\ms{Bin}\text{-}\ms{Vor}$&5.8&4.6&84.3&16.0&12.4&99.2\\
$\ms{St}\text{-}\ms{Vor}$&4.9&4.6&80.1&14.3&10.7&98.9\\
$\ms{Fb}\text{-}\ms{Vor}$&75.4&72.4&4.0&55.9&63.0&12.5\\
$\ms{Bin}\text{-}\ms{Lag}$&3.0&1.5&56.3&5.1&4.6&92.3\\
$\ms{St}\text{-}\ms{Lag}$&3.0&1.6&61.4&5.9&5.5&93.9\\
$\ms{Fb}\text{-}\ms{Lag}$&99.0&99.1&27.7&92.5&95.9&5.3
                 \end{tabular}\,
                 
                 \begin{tabular}{l|p{29pt}p{29pt}p{29pt}p{29pt}p{29pt}p{29pt}}
			 \multicolumn{7}{c}{$\TP$}\\
			 $\ms H_0\backslash \ms H_1$\hspace{-.15cm}            & $\ms {Bin}\text{-}\ms{Vor}$& $\ms{St}\text{-}\ms{Vor}$& $\ms{Fb}\text{-}\ms{Vor}$& $\ms {Bin}\text{-}\ms{Lag}$ & $\ms{St}\text{-}\ms{Lag}$ &$\ms{Fb}\text{-}\ms{Lag}$ \\
                         \hline

$\ms{Bin}\text{-}\ms{Vor}$&6.2&91.7&100.0&9.7&85.8&100.0\\
$\ms{St}\text{-}\ms{Vor}$&92.0&3.2&100.0&95.9&6.5&100.0\\
$\ms{Fb}\text{-}\ms{Vor}$&100.0&100.0&5.6&100.0&100.0&25.3\\
$\ms{Bin}\text{-}\ms{Lag}$&4.2&93.2&100.0&3.7&87.9&100.0\\
$\ms{St}\text{-}\ms{Lag}$&76.0&4.4&100.0&84.9&5.0&100.0\\
$\ms{Fb}\text{-}\ms{Lag}$&100.0&100.0&19.4&100.0&100.0&5.2
                 \end{tabular}
         \end{center}
 \end{table}

 That is, for each test statistic $T$ and null model, we compute the empirical mean and variance from these realizations. Assuming a normal distribution, we can construct asymptotic confidence intervals that yield the acceptance region of the tests. For evaluating the powers of the test, we generate a new set of 1,000 realizations per model. The hypothesis that one of these realizations is generated by a tessellation model $M_0$ is rejected, if the test statistic $T$ for the data falls outside the confidence interval of $T$ under model $M_0$. The power of the test is then obtained as the proportion of simulations where the null hypothesis is rejected.

%APROX NORM
Table \ref{pow_mark_tab} summarizes these findings. 
First, we observe that the rejection rates along the diagonal are mostly fairly close to the nominal 5\%-level. 
Second, although the face areas are a relatively elementary characteristic of the data, we see that the test $\TA$ is already quite powerful in a variety of testing scenarios. This is most clearly seen when distinguishing between the models $\ms {Bin}$-$\ms{Vor}$ and $\ms{St}$-$\ms{Lag}$, and between the models $\ms{St}$-$\ms{Vor}$ and $\ms{St}$-$\ms{Lag}$. Although, in many cases the inradii- and persistence-based statistics $\TF$ and $\TP$ lead to lower rejection rates, we see a substantial improvement for the cases where $\ms{Bin}$-$\ms{Lag}$ and $\ms{St}$-$\ms{Lag}$  are the null models and $\ms{Bin}$-$\ms{Vor}$ and $\ms{St}$-$\ms{Vor}$ are the alternatives, respectively. This already gives an indication that these new statistics can provide additional insights in settings where the traditional approaches fail. In the real data example in Section \ref{foam_sec} below, we will see an additional example for the added value of the persistence-based characteristics.

\section{Application to foam data}
\label{foam_sec}

%DAT SET
We now apply the testing methodology to a dataset from materials science describing the structure of an open aluminium alloy foam. A cubic sample of edge length 40\,mm is spatially imaged by micro computed tomography at voxel size $29.44\,\mu$m. The foam cells show an elongation in $y$-direction \cite{AWM}. Hence, we consider a rescaled version of the image to obtain an isotropic cell system. Cells are reconstructed by using the watershed transform, see \cite{AWM} for details. The resulting dataset consists of 962 cells in a $[1,657]\times [1,533] \times [1, 642]$ voxel sampling window. To represent this dataset in the form of a tessellation, we compute a Laguerre approximation by the method from \cite{LagRec}. Figure \ref{secf_fig} shows a 2D section of the original dataset and the corresponding section in the Laguerre tessellation. This comparison illustrates that the Laguerre approximation gives a good representation of the original dataset.

    Using the test statistics presented in Section \ref{sim_sec}, we want to perform goodness-of-fit tests for several tessellation models for the foam structure.    More precisely, we will investigate the face areas, face inradii as well as the persistences based on the tessellation-adapted filtration. In order to avoid edge effects, we restrict our analysis to the 356 cells that are entirely contained in the sampling window. 

 \bef[!ht]
    \includegraphics[angle=90, width=.45\textwidth]{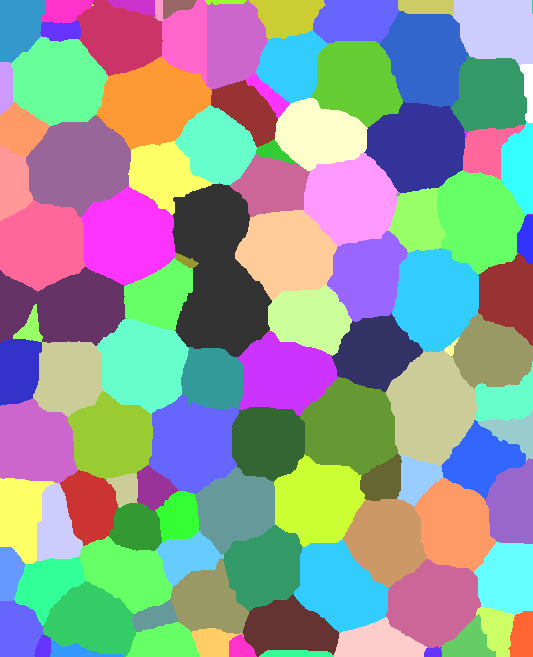}\;\;
    \includegraphics[angle=90, width=.45\textwidth]{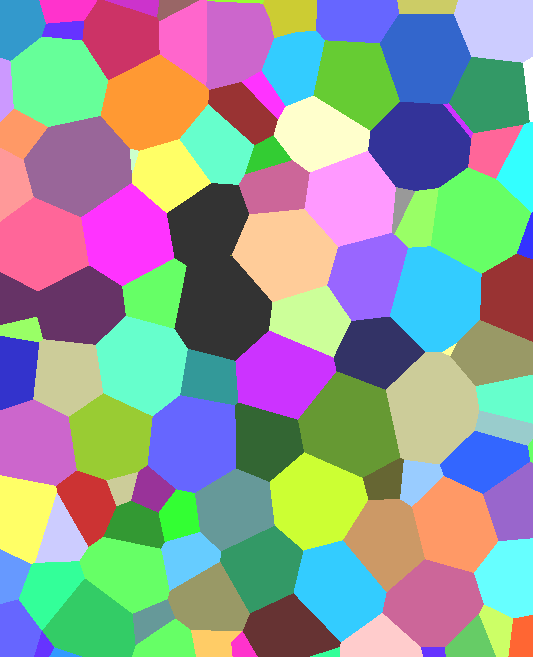}
	\caption{2D section of the original foam data (left) and the fitted Laguerre model (right).}
    \label{secf_fig}
    \enf
   
    We compare characteristics of the foam dataset to those of four different stochastic models of increasing complexity.   In the first two models, we consider  two standard tessellation models from stochastic geometry, namely Binomial- and Strauss-Voronoi tessellations.    As a prototypical example, we consider a highly repulsive point process with $\g = 0.01$ and $R = 60$. We first construct the stochastic models by using periodic boundary conditions, and then only retain the cells contained entirely in the sampling window. We choose the intensities of tessellation models so as to give in expectation 356 remaining cells. As more refined models, we consider a Laguerre tessellation based on a forced-biased sphere packing and its relaxation by the Surface Evolver \cite{Brakke1992}. Both models are fitted to first and second moments of geometric characteristics of the foam cells, see \cite{AWM} for details.

%
    %AREA
    %https://www.overleaf.com/project/5fb63ecf2f82404f0621287a
    \subsection{Exploratory analysis}
    \label{expld_sec}
    In a first step, in Figure \ref{areas_fig} (left), we compare   the face-area distributions of the model and of the data. We observe a clear difference between the two point-process based Voronoi models and the data. Both the Binomial-Voronoi and the Strauss-Voronoi tessellations exhibit a high variance: a substantial proportion of faces may exhibit a very small or a very large area. The main differences between the dataset and the force biased model are that the model still has too many faces with very small areas, and is also too sharply concentrated around the mode. We also see that in the current example, the relaxation through the Surface Evolver is slightly too aggressive in the sense that small face areas are removed entirely in this configuration, whereas there are still some faces with small areas in the dataset.

    A similar effect is also visible for the face inradii as shown in Figure \ref{areas_fig} (middle). Here, the Binomial-Voronoi and Strauss-Voronoi models exhibit a larger variability of the face inradii than the foam dataset. The sphere packing and surface evolver models lead to inradii distributions that are far more similar to the one observed in the data. However, in both models the concentration around the mode is a bit more pronounced than what we observe in the data.
    
   Finally, in Figure \ref{areas_fig} (right), we consider the total persistences. Again, we see that the Voronoi and Strauss models produce substantially too many small values. In contrast, the surface evolver model essentially removes all small persistences, whereas there are still some in the foam dataset. When considering just the persistences, the sphere packing model results in a good approximation to the foam data.

\bef[!ht]
    \includegraphics[width=.99\textwidth]{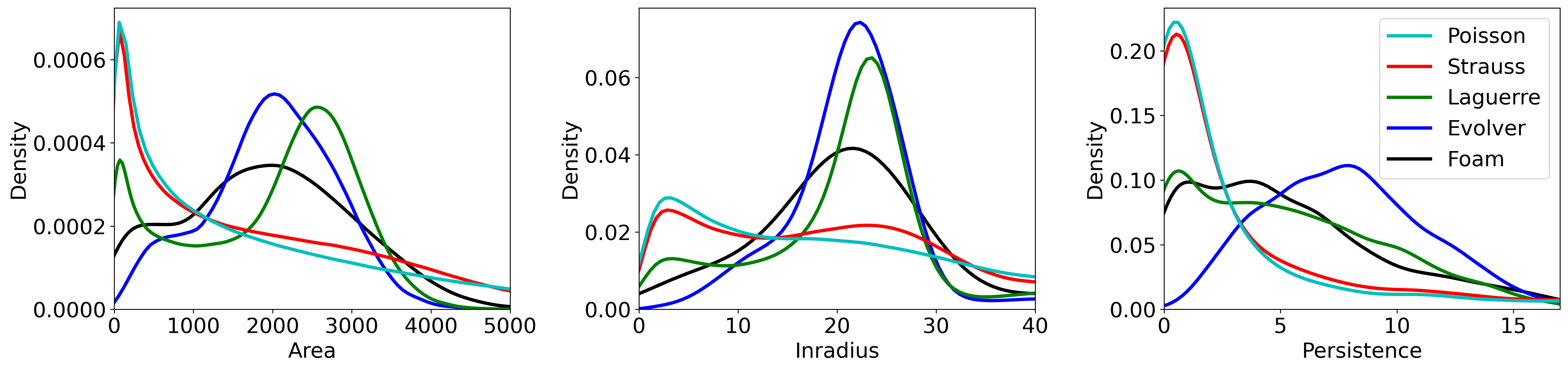}
	\caption{Kernel-density estimates for the face areas (left), face inradii (middle) and total persistences (right).}
    \label{areas_fig}
    \enf

%
%VOR
%
%\subsubsection{Persistence diagram}
Finally, we proceed to the persistence diagram on the tessellation-adapted filtration. In comparison to the face inradii, now a slightly more sophisticated analysis is needed, since the features with respect to this filtration have both non-zero birth- and death times. Nevertheless, Figure \ref{vert_fig} reinforces the previous observations: whereas the data is strongly concentrated, both the Voronoi and the Strauss model exhibit substantial fluctuations and noise. We also see that the evolver model produces fewer small persistences than observed in the foam dataset.
	\bef[!ht]
	    \centering
	    \includegraphics[width = .99\textwidth]{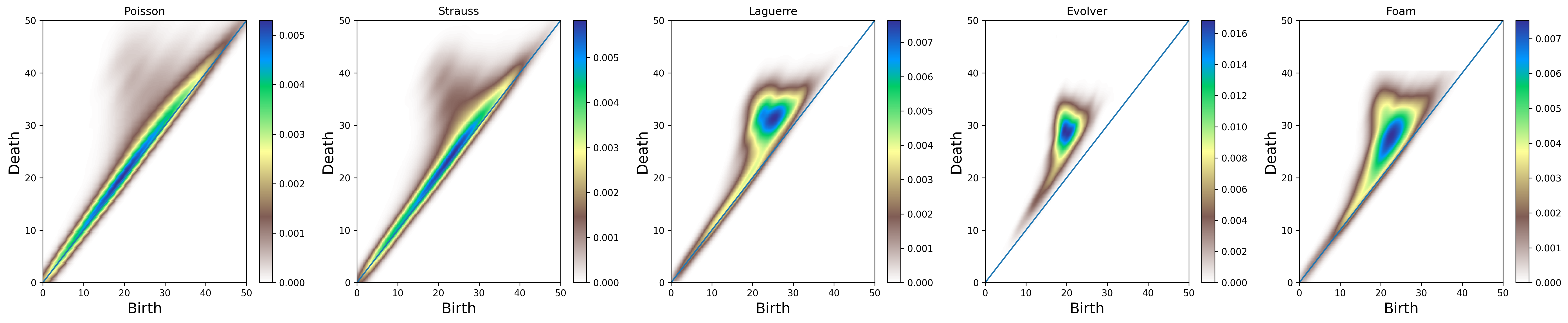}
		\caption{Persistence diagrams for the tessellation-adapted filtration.}
	    \label{vert_fig}
	    \enf

	    %
	    %FIT SEC
	    %
	    \subsection{Goodness of fit tests}

 Finally, we apply the goodness-of-fit tests from Section \ref{norm_sec} to the present datasets. To that end, we are guided by the exploratory analysis from Section \ref{expld_sec}, and choose the thresholds  $a_{\ms{Area}}, a_{\ms{I}}, a_{\ms{Pers}}$ to be the $40\%$, $60\%$, and $70\%$ quantiles of the test statistics $\TA$, $\TF$ and $\TP$ under the null model. To compute the mean and variance under the Binomial and Strauss models we generated 1,000 realizations. Due to the high computational costs, for the sphere packings we used only 100 realizations. 
 
In Table \ref{pow_dat_tab}, we present the $z$-scores for the goodness of fit test of the considered null models to the specific dataset. First, we see that the null-models of a Binomial-Voronoi or a Strauss-Voronoi tessellation are clearly rejected for all tests. Moreover, we stress that the highest $z$-scores occur for the persistence-based test statistic. This hints at the potential of TDA-based methods to assess the fit of classical random spatial tessellation models to data from materials science. Moreover, we see that the more advanced sphere-packing and evolver models lead to much lower $z$-scores, indicating by relying on such models it is possible to obtain substantially more realistic representations of the data set in contrast to the traditional standard models from stochastic geometry. As observed in Section \ref{expld_sec} for the specific example, the removal of the small surfaces in the surface evolver is too aggressive, thereby leading to higher rejection rates than for the standard sphere packing.
	   \newcolumntype{P}[1]{>{\centering\arraybackslash}p{#1}}
	 \begin{table}[!htpb]
		 \begin{center}
			 \caption{$z$-scores ($p$-values) for the test statistics $\TA$, $\TF$ and $\TP$ computed on the null models and the dataset. The rows correspond to the null models specified at the beginning of Section \ref{foam_sec}.}
                 \label{pow_dat_tab}
                 \begin{tabular}{l|P{99pt}P{99pt}P{99pt}P{99pt}}
			 $\ms H_0\backslash$Statistic\hspace{-.1cm}            & $\TA$& $\TF$ & $\TP$ \\
                         \hline\\[-2ex]
			 Binomial &12.8 $(1.6\times 10^{-37})$&14.94 $(1.9\times 10^{-50})$&35.94 $(7.5\times 10^{-283})$\\
			 Strauss &11.84 $(2.4\times10^{-32})$&14.36 $(8.9\times 10^{-47})$&32.48 $(2.4\times 10^{-231})$\\
			 Laguerre &6.92 $(4.5\times10^{-12})$&1.85 $(0.065)$&1.33 $(0.18)$\\
			 Evolver &12.42 $(2.1\times10^{-35})$&4.03 $(5.5\times 10^{-5})$&20.47 $(3.6\times 10^{-93})$
                 \end{tabular}\;
         \end{center}
 \end{table}

\section{Concluding remarks}
\label{sec:cr}

In this article, we have developed a statistical methodology to test the goodness of fit of Gibbs-based Voronoi and Laguerre tessellations based on statistics derived from a tessellation-adapted form of persistent homology. We analyzed the validity and performance of such tests for finite sample sizes in a simulation study. Finally, we illustrated how to apply our methods to a spatial tessellation originating from real foam data.

Our simulation study strongly suggests that the asymptotic normality of the test statistics should hold under substantially less restrictive assumptions than the ones imposed at the moment. In particular, it would be attractive to treat Gibbs processes with a fixed number of generators, to loosen the intensity constraint, and to allow for more flexible boundary conditions as well as possibly unbounded Laguerre marks. Furthermore, we believe that the constraint to the approximate $M$-bounded and $M$-localized Betti numbers is an artifact of our proofs and should not be needed from a conceptual perspective. Moreover, we believe that many of the methods presented here can be generalized to other tessellation models such as Johnson-Mehl tessellations.

Concerning the simulation and data study, we have so far not attempted to search for the most powerful test statistic. We have also not yet developed a clear heuristic, which test statistic is the most suitable for a given dataset. Moreover, for the dataset studied in Section \ref{foam_sec}, all of the test statistics clearly reject inappropriate null models. It would be interesting, to consider also datasets where the topological discrepancies from the null model are more subtle so that the refined persistence-base methods can unfold their full potential.

\section*{Acknowledgements}
We thank Anne Jung (Saarland University) for providing the foam sample and Christian Jung (University of Kaiserslautern) for computing the Laguerre approximation.

\bibliography{lit}
\bibliographystyle{abbrv}
\newpage

\setcounter{page}{1}
\section*{Supplementary material}
\section{Proof of Theorem \ref{clt}}
\label{clt_sec}

The proof idea is to apply the general CLT for Gibbs point processes, \cite[Theorem 2.3]{gibbs}.  To that end, we need to implement three steps: i) express $\b_{ n}^{\ms e, M, s}$ and $\b_{p - 2, n}^{M, b, d}$ as a sum of score functions, ii) verify the exponential stabilization condition \cite[Condition (2.4)]{gibbs}, and iii) verify the moment condition \cite[Condition (2.7)]{gibbs}. 

In order to approach these tasks, we note that \cite[Section 5.2(ii)]{gibbs} treats the case of edge lengths in Voronoi tessellations. Moreover, \cite[Theorems 2.4, 2.5]{flimmel} deals explicitly with the case of the surface area and the distribution function of the cell volume. In addition, it is observed that the techniques extend to a variety of further score functions.  Hence, the main step in the proof of Theorem \ref{clt} is to adapt and combine the arguments so that they apply to $\b_n^{\ms e, M, s}$ and $\b_{p - 2, n}^{M, b, d}$. Although it will lead to some redundancy with \cite{flimmel,gibbs}, we present here some details of this program to make the overall presentation more accessible. Now,  we first reproduce the exponential stabilization and moment conditions, and then verify them for $\b_n^{\ms e, M, s}$ and $\b_{p - 2, n}^{M, b, d}$. We present the moment condition by relying on a fixed nonnegative score function $\xi(x, \AA)$, capturing some form of interaction of a point $x \in \Rm$ with a finite configuration $\AA \su \Rm$. 
\medskip

%
%MOMENT
%
\ni{\bf $k$-moment condition \cite[Condition (2.7)]{gibbs}.}
Fix $k \ge 0$ and a Hamiltonian $\Psi$ of class $\bs\Psi$. Then, the \emph{$k$-moment condition is satisfied} if 
$$\sup_{n \ge 1} \sup_{\cx \in \wc Q_n}\sup_{ \AA \su \Rm \text{finite}} \E\big[\xi\big(\cx, \XX_n \cup  \AA \big)^k\big] < \ff$$
holds for all  $\t < |B_1(o)|^{-1}\exp(\a m_0^\Psi)(r_0 + 1)^{-p},$ where $\wc Q_n(x) := (x + [-n/2, n/2]^p) \times [0, \mu]$.
\bigskip

To formulate the condition on exponential stabilization, we also need the concept of a Poisson-like processes. We write $B_s(x) := \{y \in \R^p \co |y - x| \le s\}$ for the Euclidean ball of radius $s > 0$ centered at $x \in \R^p$ and put  $\Bm_s(x) := B_s(x) \times [0, \mu]$.

%
%POISSON LIKE
%
\ni{\bf Poisson-like processes \cite[Condition (2.1)]{gibbs}.} Consider the marked Poisson process $\PP_{\t, \Q}$ as introduced in Section \ref{cond_sec}.  Then, a point process $\YY$ is called  \emph{Poisson-like} if the following two conditions are satisfied:
\been
\im 
$\P(\YY(B)\ge n ) \le \P(\PP_{\t, \Q}(B) \ge n)$
holds for every $n \ge 0$ and Borel set $B \su \Rm$;
\im there exist $c_{\ms{PL}}, r_{\ms{PL}} > 0$ such that 
$$\P\big(\YY \cap \Bm_s(x) = \es \ba \YY \sm \Bm_s(x) = \mc E_s(x)\big) \le \exp(-c_{\ms{PL}}s^p)$$
holds for every $r \ge r_{\ms{PL}}$, $x\in\R^p$ and locally finite $\mc E_s(x) \su \Rm \sm \Bm_s(x)$.
\enen
By \cite[Lemma 3.3]{gibbs}, Gibbs point processes with Hamiltonian in $\bs \Psi^*$ and that satisfy $\t < |B_1(o)|^{-1}\exp(\b m_0^\Psi)(r_0 + 1)^{-p}$ are Poisson-like.

\bigskip
%
%EXP STAB
%
\ni {\bf Exponential stabilization \cite[Condition (2.4)]{gibbs}.} Let $\YY$ be a simple point process and for $\cz\su\Rm$ write $\YY^{\cz}:= \YY \cup \cz$. Then, the definition of exponential stabilization for the score function $\xi$ relies on the notion of a \emph{radius of stabilization} $R = R\big(\cx, \YY^{\cz}\big)$, which is determined by $\cx \in \Rm$, $\cz\su \Rm$ of cardinality at most 1, and the realization $\XX$, and which satisfies that for almost all realizations of $\YY$ that 
$$\xi\big(\cx, \YY^{\cz} \cap \Bm_R(x)\big) = \xi\big(\cx, (\YY^{\cz} \cap \Bm_R(x)) \cup \AA\big)$$
for every finite $\AA \su \Rm \sm \Bm_R(x)$.  Having defined the radius of stabilization, we next introduce the conditional tail probability
$$t(s, \e) := \sup_{y \in \R^p}\sup_{\substack{\cz\su \Rm \\ \#\cz \le 1}} \P\big(\sup_{\cx \in \YY^{\cz}\cap  \Bm_\e(y)} R(\cx, \YY^{\cz} )> s \big).$$
Then, we demand that $\limsup_{\e\to 0}\limsup_{s \to\ff} s^{-1}{\log t(s, \e)} < 0$ hold for every Poisson-like $\YY$.

%
%EBB
%
\subsection{Edge-based persistent Betti numbers}
To represent $\b_n^{\ms e, M, s}$ as a sum of score functions, we first associate each bounded face $f\in \Xi_n^{(2)}$ to precisely one of the 3D Laguerre cells  $C(\cX_i, \XX_n)$ incident to $f$. To be more precise, among all centers $X_i$ whose associated cell $C(\cX_i, \XX_n)$ is incident to $f$, we let $X(f, \XX_n)$ denote the lexicographic minimum.  This correspondence makes it possible to express $\b_n^{\ms e, M, s}$ as a sum of scores associated with the points of $\XX_n$. More precisely, 
$$\b_n^{\ms e, M, s} = \sum_{\cX_i \in \XX_n} \b^{\ms e, M, s}\big(\cX_i, \XX_n\big),$$
where 
$$\b^{\ms e, M, s}(\cX_i, \XX_n):= \#\Big\{f \in \Xi^{(2)}_{n, M}\co \rfc(f) > s \text{ and } X_i = X(f, \XX_n)\Big\}.$$

%
%MOM PRF
%
We begin with a geometric auxiliary result ensuring the boundedness of Laguerre cells. 

%
%BCELL LEM
%
\bel[Boundedness of cells]
\label{stab_lem}
Let $\AA\su \Rm$ be locally finite and let $a \ge \mu$. Assume that $\Qm_a(az) \cap \AA \ne \es$ for every $z \in \Z^p$ with $|z|_\ff =  2p + 2$. Then, 
$$C(\cx, \AA) \cap \big(Q_{(4p + 3)a}(o)\sm Q_{(4p + 1)a}(o)\big) = \es$$
for every  $\cx \in \AA \cap \big(\Qm_a(o) \cup(\Rm \sm \Qm_{(8p  + 5)a}(o))\big)$.
%\su Q_{(4p + 3)a}(o)$ for every $\cx  \in \AA \cap \Qm_a(o)$.
\enl
\bep
Let $P \in Q_a(az)$ for some $z \in \Z^p$ with $|z|_\ff = 2p + 2$. Then, for $x' \in Q_a(az)$ we have
$$|P - x|^2 - \mu^2 \ge 4p^2a^2 - \mu^2 > pa^2 \ge |P - x'|^2,$$
thereby showing  $P \not \in C((x, r), \AA)$,   as asserted.
\enp

Now, we elucidate on how to verify exponential stabilization and the $k$-moment condition.

\bep[Proof of exponential stabilization]
For $a > 0$ choose $z_0 := z_0(a) \in \Z^p$ such that $x \in Q_a(az_0)$. Furthermore, let 
$$E_{a, 0} := E_{a, 0}(x) :=\Big\{Q_a(az) \cap \YY \ne \es \text{ for all $z\in\Z^p$ with $|z  - z_0|_\ff = 2p + 2$}\Big\}$$
denote the event that all boxes in the $(2p + 2)a$ neighborhood of $Q_a(az_0)$ contain at least one point from $\YY$. Then, Lemma \ref{stab_lem} implies that $C(\cx, \YY^{\cz}) \su Q_{(4p + 1)a}(az_0)$ for every $\cx \in \YY^{\cz} \cap \Bm(y)$. Moreover, again by Lemma \ref{stab_lem}, changing $\YY^{\cz}$ outside $\Qm_{(8p + 5)a}(az_0)$ does not influence $C(\cx, \YY^{\cz})$. Hence, under the event $E_{a, 0}$ the stabilization radius is at most $(8p + 5)\sqrt p a$. We conclude the proof by noting that since $\YY$ is Poisson-like the probability of the complement of $E_{a, 0}$ decays at exponential speed in $a^p$.
\enp

\bep[Proof of the $k$-moment condition]
To verify the moment condition we show that $\b^{\ms e, M, s}\big(\cx, \XX_n \cup \AA \big) \le 16\pi M^2/s^2$. To prove this assertion, note that by definition of $\ms{ecc}(f)$, we only need to consider faces with inradius at least $s$ that are contained in the ball $B_{M}(x)$. In particular, the surface area of each such face is at least $\pi s^2$. Hence, since $C\big(\cx, \XX_n \cup \AA\big)$ is convex, 
$$\b^{\ms e, M, s}\big(\cx, \XX_n \cup \AA \big) \le 4\pi (2M)^2/(\pi s^2)= 16\pi M^2/s^2, $$
as asserted.
\enp

%
%MLB
%
\subsection{$M$-localized persistent Betti numbers}
The arguments for the $M$-localized persistent Betti numbers are similar to those for the edge-based Betti numbers. Hence, to avoid redundancy we will only focus on the parts, where substantial differences occur. First, we represent $\bbnd$ as a sum of score functions, by setting
$$\bbnd = \sum_{\cX_i \in \XX_n} \bbd\big(\cX_i, \XX_n\big),$$
where 
$$\bbd(\cX_i, \XX_n):= \#\big\{f \in \Xi_{n, M}^{',(p - 1)}\co b(f) \le b,\, d(f) \ge d\, \text{ and } X_i = X(f, \XX_n)\big\}.$$

%
%MOM PRF
%
Next, we elucidate on how to verify exponential stabilization and the $k$-moment condition. 

\bep[Proof of exponential stabilization]
In contrast to the edge-based case we need to stabilize the cells contained in an $M$-neighborhood of $C(\cx, \YY^{\cz})$. To achieve this goal, we assume from now on that $a \ge 2M$. Again, we fix $z_0 := z_0(a) \in \Z^p$ such that $y \in Q_a(az_0)$. Furthermore, set 
$$E_{a, 1} := E_{a, 1}(x) :=\Big\{Q_a(az) \cap \YY \ne \es \text{ for all $z\in\Z^p$ with $|z - z_0|_\ff \le 4p + 6$}\Big\}.$$
Now, we may again apply Lemma \ref{stab_lem} to deduce that changing $\YY^{cz}$ outside $\Qm_{(12p + 16)a}(az_0)$ does not influence $C(\cx, \YY^{\cz})$ or one of the cells in an $M$-neighborhood. Again, since the Poisson-likeness of $\YY$ implies that the complements of the events $E_{a, 1}$ decay at exponential speed in $a^p$, we conclude the proof.
\enp

\bep[Proof of the $k$-moment condition]
As in the edge-based case, we actually derive a deterministic upper bound on $\bbnd$. Indeed, relying again on the definition of the eccentricity, we deduce that the relevant faces are contained in the ball $B_M(x)$. Hence, using that by definition of $\Xi_{n, M}^{',(p  - 1)}$, we have a lower bound on the inradii, we conclude as in the edge-based case.
\enp

\section{Perfect simulation of Gibbs processes and concepts of stabilization}
\label{gibbs_sec}

To prepare the proof of Theorem \ref{fclt},  it is instructive to review the perfect simulation of Gibbs processes, which originated in \cite{garcia}. To make our arguments self-contained, we summarize the main steps from the excellent exposition in \cite[Section 3.2]{gibbs}.

%
%BOUND DOM
%
Let $\{\r_n(t)\}_{t \in \R} := \{\r_{W_n}(t)\}_{t \in \R}$ denote a stationary and homogeneous free birth-and-death process in the marked window $\Wm_n := W_n \times [0, \mu]$. That is, particles have rate-1 exponentially distributed lifetimes and are born according to the space-time intensity $\t \d x\d t$ for some $\t > 0$.

%
%TRIM
%
This free process gives rise to a trimmed process $(\g_n(t))_{t \in \R}$ as follows. Let $T_-(t) := \sup\{s \le t \co \r_n(s) =\es\}$ be the last time before $t$, where the free birth process is empty, and define $\g_n(T_-(t)-) := \es$. Then, for $T_-(t) \le s \le t$ the trimmed process evolves as follows. A birth at time $s \in \R$ and marked location $\cx \in \Wm_n$ in the free process is accepted with probability $\exp\big(-\a \Delta (\cx, \g_n(s-))\big)$. In the case of rejection, no particle is born at time $s$.

%
%ANC
%
We now describe an alternative construction to encode the near-independence in spatially distant regions.
Suppose in the free process a birth occurs at location $x \in W_n$ and time $t \in \R$. Then, the ancestors of this birth consist of all births within the interaction range $B_{r_0}(x)$ that are still alive at time $t$. We can then form the ancestors of the ancestors, and by repeating this construction arrive at the \emph{ancestor clan}.  We let $\A_{D, n}(t) \su W_n$ denote the ancestor clans for all births within a set $D \su W_n$ that are alive at time $t$ and put $\A_{D, n} := \A_{D, n}(0)$. Our proof will rely on two key findings from the analysis in \cite[Section 3.2]{gibbs}.

%
%INDEP
%
 First, the ancestor clans help to encode independence between distant spatial regions. More precisely, for $D_1\su D_1' \su W_n$ and $D_2 \su  D_2' \su W_n$ with $D_1' \cap D_2' = \es$ and nonnegative measurable $f,f'$, the random variables $f(\XX \cap D_1) \one\{\A_{D_1, n} \su D_1'\}$ and $f(\XX \cap D_2) \one\{\A_{D_2, n} \su D_2'\}$ are independent. Here, to ease notation, we put $\XX\cap D = \big\{(X_i, R_i) \in \XX \co X_i \in D\big\}$.

%
%EXP
%
Second, it is highly unlikely that ancestor clans are far larger than the domain from which they were constructed.  More precisely, by \cite[Equation (3.6)]{gibbs}, there exists a constant $\ca > 0$ such that for all $v \ge 1$ 
\begin{align}
	\label{anc_bound_eq}
	\sup_{a, s \ge 1}\sup_{n \ge 1}\max_{z \in \on}a^{-p}\P\big(\diam(\A_{W_s(z), n} \ge v + \sqrt p a\big)  \le \ca\exp(-v/\ca),
\end{align}
where $W_s(z) := z + [-\tf s2 + \tf12, \tf s2 + \tf12]^p$. Finally,  we introduce the \emph{ancestral stabilization radius} 
%STAB
$$\RAi_n(A):= \min\big\{m \ge 1\co \bigcup_{z\in A}\A_{W_1(z), n} \su W_m(A)\big\}$$
of $A \su \on$ as the smallest $m \ge 1$ such that the ancestor clans $\A_{A, n}$ are contained in  $W_m(A) := \bigcup_{z\in A}W_m(z)$. To combine the ancestor clans with the stabilization radii from the tessellation on the resulting Gibbs process, we define the \emph{combined internal stabilization radius}
$$\Rci_n(A) := \RAi_n(W_{R_n(A)}(A)),$$
where we set $R_n(A) := \sup_{\cx \in W(A) \times [0, \mu]}R(\cx, \XX_n)$ with $W(A) := W_1(A)$. It is also helpful to be aware of dependencies in the converse direction: to what extent do changes in the free birth-and-death process inside $W(A)$ influence the tessellation at large distances? In order to capture this dependence, we introduce the \emph{total external stabilization radius} along the same lines as the original one. More precisely, we let
$$\r_{m, \AA, \BB, n} :=  \AA \cup \BB \cup\bi(\r_n \cap (W_m(A) \sm W(A))\bi)$$
be the birth-death process obtained from $\r_n$ by placing a finite birth-and-death process $\AA$ outside of $W_m(A)$ and a finite birth-and-death process $\BB$ inside $W(A)$. Then, $\Xi_n$ is \emph{externally stabilized} within $A$ at distance $m$ if
$$C_n(\cX_i, \XX(\r_{m, \AA, \es, n}))  = C_n(\cX_i, \XX(\r_{m, \AA, \BB, n})) $$
holds for every finite birth-and-death processes $\AA$ outside $W_m(A)$, $\BB$ inside  $W(A)$ and every center point $\cX_i$ outside $\Wm_m(A)$. Finally, we introduce the \emph{total external stabilization radius} $\Rce_n(A)$ as the smallest integer $m$ such that $\Xi_n$ is externally stabilized within $A$ at distance $m$. Using \eqref{anc_bound_eq}, we note that under the exponential stabilization as described in Section \ref{clt_sec} the tails of both $\Rce_n$ and $\Rci_n$ decay at a stretched exponential speed.

\section{Proof of Theorem \ref{fclt}}
\label{fclt_sec}
Henceforth, omit the feature dimension $p - 1$ from the notation. Proving a functional CLT requires two steps: i) asymptotic normality of the marginals, and ii) tightness as a stochastic process. As Theorem \ref{clt} deals with the normality of the marginals, it remains to establish tightness. To that end, we rely on the Chentsov-type moment criterion from \cite{bickel}. We first explain more precisely how to proceed for the $M$-localized persistent Betti number since this is the more complex case. For a block $E = [\tbm, \tbp] \ti [\tdm, \tdp]$, we let
$$\b_n^M(E) := \b_n^{M, \tbp, \tdp} - \b_n^{M, \tbm, \tdp} - \b_n^{M, \tbp, \tdm} + \b_n^{M, \tbm, \tdm}$$
denote the increment of the $M$-bounded persistent Betti number in that block. One potential issue is that the expression $\b_n^{M, \tbp, \tdm}$ is ill-defined for $\tbp > \tdm$. However, we can use this situation to our advantage and put in that case
$$\b_n^{M, \tbp, \tdm} := \b_n^{M, \tbp, \tbp} - \b_n^{M, \tdm, \tbp} + \b_n^{M, \tdm, \tdm}.$$
Thus, $\b_n^M([\tdm, \tbp]^2) = 0$ so that we may henceforth assume that $\tdm \le \tbp$.

Then, by \cite[Theorem 3]{bickel}, the goal is to prove that 
\begin{align}
	\label{lav_eq}
	\sup_{n \ge 1}\sup_{E \su \ot}\f{\E[\bar\b_n^M(E)^4]}{n^{2p}|E|^{9/8}} < \ff, 
\end{align}
where $\bar\b_n^M(E) := \b_n^M(E) - \E[\b_n^M(E)]$ denotes the centered increment, and where here and in the following we tacitly assume that the supremum is taken over blocks $E \su \ot$. For the edge-based persistent Betti numbers, we proceed in the same manner except that now blocks are just 1-dimensional, i.e., $E \su [0, S]$. Hence, it suffices to verify the simpler one-parameter condition from  \cite[Theorem 13.5]{billingsley}. 

%GRID
To verify the Chentsov condition, we will build on the blueprint developed in \cite{cyl} for the setting of networks in cylindrically growing domains. That is, we proceed in two steps: \emph{reduction to a grid} and \emph{variance and cumulant  bounds}.
First, we rely on a trick from \cite{davydov} in order to reduce the verification to blocks that are not too small. We state this reduction step now and defer its proof to Section \ref{grid_sec}. Henceforth, we say that a block $E \su [0, S]$ is \emph{$n$-big} if $|E| \ge n^{-9/4}$; a block $E \su [0, S]^2$ is \emph{$n$-big} if $|E| \ge n^{-2p}$.
\vspace{.9cm}

%
%GRID PROP
%
\bepr[Reduction to $n$-big blocks]
\label{grid_prop}
\phantom{}
\been
\im If condition \eqref{lav_eq} holds for all $n \ge 1$ and $n$-big blocks $E \su [0, S]$, then the processes $\{\b_n^{\ms e, r}\}_{r \le S}$ are tight in the Skorokhod topology.
\im If condition \eqref{lav_eq} holds for all $n \ge 1$ and $n$-big blocks $E \su \ot$, then the processes $\{\b_n^{M, b, d}\}_{b, d \le S}$ are tight in the Skorokhod topology.
\enen
\enpr

The second key ingredient are bounds on the variance $\Var(\b_n(E))$ and the fourth-order cumulant $c^4(\b_n(E))$.   
Again, we state these bounds now and prove them later in Section \ref{var_sec}.

%
%VAR BOUND
%
\bepr[Variance and cumulant bound]
\label{var_prop}
\phantom{}
\been
\im
It holds that
$$\sup_{n \ge 1}\sup_{\substack{E \su [0, S]\\ E \text{ is $n$-big}}}\f{\Var\big(\b_n^{\ms e, M}(E)\big) |E|^{-5/8}+ c^4\big(\b_n^{\ms e, M}(E)\big)}{n^3 } < \ff .$$
\im
It holds that
$$\sup_{n \ge 1}\sup_{\substack{E \su \ot \\ E \text{ is $n$-big}}}\f{\Var\big(\b_n^M(E)\big)|E|^{-5/8} + c^4\big(\b_n^M(E)\big)|E|^{-5/8}}{n^p } < \ff .$$
\enen

\enpr

%
%FCLT PROOF
%
We now elucidate the proof of Theorem \ref{fclt}.
\bep[Proof of Theorem \ref{fclt}]
Starting the proof in the setting of the edge-based persistent Betti numbers, we invoke the cumulant identity $\E[\bar\b_n^{\ms e, M}(E)^4] = 3\Var(\b_n^{\ms e, M}(E))^2 + c^4(\b_n^{\ms e, M}(E))$ and then apply part (1) of Proposition \ref{var_prop} to arrive at 
\begin{align*}
	\E[\bar\b_n^{\ms e, M}(E)^4] &= 3\Var(\b_n^{\ms e, M}(E))^2 + c^4(\b_n^{\ms e, M}(E)) \le 3c^2n^{6} |E|^{5/4} + c n^3.
\end{align*}
Now, the $n$-bigness of $E$ gives that $n^3 \le n^{6}|E|^{4/3}$, which concludes the proof. Similarly, for the $M$-localized persistent Betti numbers, 
\begin{align*}
	\E[\bar\b_n^{M}(E)^4] &= 3\Var(\b_n^{M}(E))^2 + c^4(\b_n^{M}(E)) \le 3c^2n^{2p} |E|^{5/4} + c n^p|E|^{5/8}.
\end{align*}
Now, the $n$-bigness of $E$ gives that $n^p|E|^{5/8} \le n^{2p}|E|^{9/8}$, which concludes the proof.
\enp

\subsection{Proof of Proposition \ref{grid_prop}}
\label{grid_sec}

The proof of Proposition \ref{grid_prop} consists of two key ingredients. First, as in \cite[Proposition 5]{cyl}, we leverage \cite[Corollary 2]{davydov}, thereby reducing the tightness proof of $\b_n^{\ms e, s}$ and $\b_n^{M, b, d}$ to the setting where the latter are restricted to  grids of the form $G_n := n^{-9/4}\Z$ and $G_n := n^{-p}\Z^2$, $n \ge 1$, respectively. Second, we rely on the following Lipschitz continuity property of the filtration times.  It is at this point where we use the assumption that the model contains additional noise components.
%
%CONT LEM 
%
\bel[Lipschitz continuity of filtration times]
\label{cont_v_lem}\phantom{}\\
\been
\im Let  $\{\r_i\}_{i \ge 1}$ be iid with a Lipschitz distribution function.  Consider $k \ge 1$ edges $e_1, \dots, e_k$ in $\R^p$ that form the boundary of a convex polygon $f$ in $\R^p$. Furthermore, attach to each $e_i$ the thickness $\r_i$.  Then, the function $t \mapsto k^{-1}\P\big( \rfc(f) \le t\big)$ is Lipschitz continuous on $[0, T]$ uniformly over all $k \ge 1$ and all edges $e_1, \dots, e_k$ bounding a convex polygon.
\im  Let  $\{\eta_ i\}_{i \ge 0}$ be iid uniform in $B_{h_0}$ for some $h_0  > 0$.  Then, the function $b \mapsto \P\big( \rci\big(P_0 + \eta_0 , \dots, P_q + \eta_q\big) \le b\big)$ is Lipschitz continuous on $[0, T]$ uniformly over all $P_0, \dots, P_q \in  \R^p$.
\enen
\enl
%
%CONT LEM 
%
\bep
For $A, A' \su \R^p$ put $A \oplus A' := \{a + a'\co a \in A,\, a' \in A'\}$.

\ni {\bf Part (1).}
	Let $0\le r_- \le r_+$ and let $f$ be a face bounded by edges $\{e_1, \dots, e_k\}$. Assume that the face $f$ is covered by $\bigcup_{i \le k} e_i \oplus B_{\r_i + \trp}(o)$ but not by $\bigcup_{i \le k} \big(e_i\oplus B_{\r_i + \trm}(o)\big)$. Hence, for some $i_0 \le k$, the face $f$ is covered by
	$$\big(\bigcup_{i \le i_0} e_i \oplus B_{\r_i + \trp}(o)\big) \cup \big( 
\bigcup_{i > i_0} e_i \oplus B_{\r_i + \trm}(o)\big)$$
	but not by
	$$\big(\bigcup_{i \le i_0 - 1} e_i \oplus B_{\r_i + \trp}(o)\big) \cup \big(\bigcup_{i \ge i_0}  e_i \oplus B_{\r_i + \trm}(o)\big).$$
	In particular, $s_{i_0}\in [\r_{i_0} + \trm, \r_{i_0} + \trp]$, where $s_{i_0}$ is the distance of the most distant point in
	$$\big(\bigcup_{i \le i_0 - 1} e_i \oplus B_{\r_i + \trp}(o)\big) \cup \big(\bigcup_{i > i_0}  e_i \oplus B_{\r_i + \trm}(o)\big).$$
	from the edge $e_{i_0}$. Therefore, writing $c_{\ms R}$ for the Lipschitz constant of $\P(\r_i \le t)$,  we arrive at the asserted
	$$\P(\rfc(f) \in [r_-, r_+]\ba \XX_n) \le \sum_{i \le k} \P\bi(\r_i \in [\trm - s_{i_0}, \trp - s_{i_0}]\bi)  \le  c_{\ms R}k (r_+ - r_-).$$
	\medskip

\ni {\bf Part (2).}
First, we express the probability in the assertion in integral form, i.e., as 
\begin{align*}
	&|B_{h_0}(o)|^{-(q + 1)} \int_{B_{h_0}(o)^{q + 1}} \one\big\{ \rci\big(P_0 + x_0 , \dots, P_q + x_ q\big) \in (b_-, b_+)\big\}\d x_0 \cdots \d x_q \\
	&\qquad\le|B_{h_0}(o)|^{-(q + 1)} \int_{B_S(o)^{q + 1}} \one\big\{ \rci\big(x_0', \dots, x_q'\big) \in (b_-, b_+)\big\}\d x_0' \cdots \d x_q'.
\end{align*}
Now, \cite[Lemma 6.10]{divol} shows the asserted upper bound $c(b_+ - b_-)$ for a suitable $c > 0$.
\enp
%
%TESS VERT LEM
%
To carry out the reduction step in Proposition \ref{grid_prop}, we need moment bounds on the number of tessellation faces in a bounded domain. More precisely, we write $\Xi^{(q)}_n(A)$ for the number of $q$-faces in $\Xi_n$ that are contained in the Borel set $A \su \R^p$. 

\bel[Moment bounds for tessellation faces]
\label{vert_lem}
Let $q\le p$. Then, there exists $c > 0$ with
\begin{align}
	\label{vert_sup_eq}
	\limsup_{s \to \ff}\sup_{n \ge 1}\sup_{z \in \on}\f{\log \P\big(\Xi^{(q)}_n(W_s(z)) > s^{2p(p + 1)}\big)}{s^c} < 0.
\end{align}
In particular, for every $K \ge 1$, 
$$\sup_{s \ge 1}\sup_{n \ge 1}\sup_{z \in \on}s^{-2p(p + 1)K}\E\big[\Xi_n^{(q)}(W_s(z))^K\big] < \ff. $$
\enl
\bep
%SUP_EQ
To verify \eqref{vert_sup_eq}, we let 
$$W_s'(z) := \big\{z \in \Z^p\co W(z) \cap W_s(z) \ne \es\big\}$$
be the discretization of $W_s(z)$. 
Now, note that $\Xi_n^{(q)}(W_s(z)) \le \XX_n\big(W_{s + R_n(W_s'(z))}(z)\big)^{p + 1}$ since cells centered outside $W_{s + R_n(W_s'(z))}(z)$ do not intersect $W_s(z)$.  Thus, 
$$\P\big(\Xi_n^{(q)}(W_s(z)) > s^{2p(p + 1)}\big) \le \P\big(\XX_n\big(W_{2s }(z)\big)> s^{2p}\big) + \P(R_n(W_s'(z)) > s).$$
To conclude the proof, note that by Poisson-likeness and exponential stabilization  both summands on the right-hand side decay at stretched exponential speed in $s$.
It remains to deduce the finiteness of the expectation from the stretched exponential decay. We may assume that $s \le n$ because $\Xi_n^{(q)}(W_s(z)) = \Xi_n^{(q)}(W_n)$ for $s > n$. Next, setting $p' := 2p(p + 1)$ gives that
\begin{align*}
	\E\big[\Xi_n^{(q)}(W_s(z))^K\big] &= K\int_0^\ff a^{K - 1}\P\big(\Xi_n^{(q)}(W_s(z)) > a\big) \d a \\
	&\le Ks^{p'K}+ K\int_{s^{p'}}^\ff a^{K - 1}\P\big(\Xi_n^{(q)}(W_s(z)) > a\big) \d a.
\end{align*}
The integral on the right-hand side remains bounded as $s \to \ff$ since $\P\big(\Xi_n^{(q)}(W_s(z)) > a\big) \le \P\big(\Xi_n^{(q)}(W_{a^{1/p'}}(z)) > a\big)$, which by \eqref{vert_sup_eq} decays at stretched exponential speed in $a$. 
\enp

%
%GRID_PROOF
%
\bep[Proof of Proposition \ref{grid_prop}]\phantom{}\\
%DAVYDOV
\ni{\bf Part (1).}	
Since the edge-based Betti numbers $\b_n^{s}$ are decreasing in $s$, the aforementioned \cite[Corollary 2]{davydov} reduces the claim to showing that 
$\E\big[ \b^{\ms e, M,   s_-}_n - \b^{\ms e, M,   s_+}_n \big]  \in o(n^{-3/2})$
for $0 \le s_+ - s_- \le n^{ -9/4}$. First, writing $C_{\ms e}$ for the Lipschitz constant from Lemma \ref{cont_v_lem}, part (1) yields that 
$$ \E\big[\b^{  s_-}_n - \b^{s_+}_n\big]\le \E\Big[\sum_{\substack{f \in  \Xi_n^{(2)}}}\one\{\rfc(f) \in [s_-, s_+]\}\Big] \le \E\Big[\sum_{\substack{f \in  \Xi_n^{(2)}}}v(f)\Big]n^{-9/4},$$
where $v(f)$ denotes the number of vertices bounding the face $f$. Since the Laguerre tessellation is normal, each vertex is contained in 4 faces, so that 
$ \E\big[\b^{\ms e, M,   s_+}_n - \b^{\ms e, M, s_-}_n\big]\le 4\E[|\Xi_n^{(0)}|]n^{-9/4}.$
Since $\E[|\Xi_n^{(0)}|]$ is of order $O(n^3)$ by Lemma \ref{vert_lem}, we obtain the claim for part (1).
\medskip

\ni{\bf Part (2).}
Again, we want to show that
$$\E\big[ \b^{M,  b_+, d}_n - \b^{M,  b_-, d}_n \big] +   \E\big[ \b^{M,  b, d_-}_n - \b^{M,  b, d_+}_n\big] \in o(n^{-p/2})$$
for $0 \le b_+ - b_-, d_+ - d_- \le n^{-p}$.  We only deal with the death times since the arguments for the birth times are analogous. 

First, by definition,
\begin{align*}
	\b^{M,  b, d_-}_n - \b^{M,  b, d_+}_n \le  \#\big\{f \in \Xi_{n}^{(q)}\co d(f) \in [d_-, d_+]\big\}.
\end{align*}
Hence, we need to bound the number of $q$-faces $f$ with circumradius $\rci(f) \in [d_-, d_+]$. Since this means that one of the $q$-simplices in $f$ has this circumradius, we deduce that
\begin{align*}
	\b^{M,  b, d_+}_n - \b^{M,  , d_-}_n \le \sum_{\cX_i \in \XX_n} \sum_{P_0, \dots, P_q \in \Xi^{(0)}_n  \cap W_M(X_i)}\one\big\{\rce(P_0 + \eta_0, \dots, P_q + \eta_q) \in [d_-, d_+]\big\}.
\end{align*}
Instead of summing over all $X_i \in \XX_n$, we first partition $W_n$ into unit boxes $W(z_1)$, \dots, $W(z_{n^p})$, and then group the points $X_i$ by the box in which they are contained. Therefore, 
$$ \b^{M,  b, d_-}_n - \b^{M, b , d_+}_n \le \sum_{\substack{i \le n^p\\P_0, \dots, P_q \in \Xi^{(0)}_n  \cap W_{M + 1}(z_i)}}\XX_n(W(z_i))\one\big\{\rci(P_0 + \xi_0, \dots, P_q + \xi_q) \in [d_-, d_+]\big\} .$$
Thus, by Lemma \ref{cont_v_lem},
$ \E\big[\b^{b, d_-}_n - \b^{b, d_+}_n \big]\le n^p(d_+ - d_-) \sup_{i \le n^p} \E\big[\XX_n(W(z_i))\Xi^{(0)}_n(W_{M }(z_i))^{q + 1}\big],$
which by Lemma \ref{vert_lem} is indeed of order $o(n^{p/2})$. 
\enp

\subsection{Proof of Proposition \ref{var_prop}}
\label{var_sec}

The key idea to establish variance and cumulant bounds is to build on the proof strategy of \cite[Theorem 3.1]{yukCLT} and express $\bar\b_n(E) $ as a sum of martingale differences. Compared to the adaptation in \cite{cyl}, the present setting is more delicate as it involves Gibbs instead of Poisson point processes. We will follow a two-fold strategy to address these challenges. First, the filtration for the martingale includes the entire free birth-and-death process instead of only the spatial point locations. Second, we invoke the exponential clustering of Gibbs processes from \cite[Lemma 3.4]{gibbs}. Let

%
%M-DECOMP
%

$$\GG_i := \s\big(\r_{W(\{z \in \on\co z  \le_{\ms{lex}} z_i\}}(t)\big)$$
denote the $\s$-algebra generated by the free birth-and-death process in the set $W(\{z \in \on\co z \le_{\ms{lex}} z_i\})$,  where for fixed $n\ge 1$, we let $z_1, z_2, \dots, z_{n^p}$ be the enumeration of the lattice points $\on$ in the lexicographic order $\le_{\ms{lex}}$.

Now, we consider the martingale-difference decomposition
\begin{align}
        \label{decomp_eq}
	\bar    \b_n^{\ms e, M}(E)  = \sum_{i \le n^p }D_{i, n}^{\ms e}(E):=  \sum_{i \le n^p }\big(\E[\b_n^{\ms e, M}(E)\ba \GG_i ] - \E[\b_n^{\ms e, M}(E)\ba \GG_{i - 1}]\big)
\end{align}
The decomposition for $\bar\b_n^M(E)$ is defined similarly. 
The first key element in the proof of Proposition \ref{var_prop} are moment  bounds for the increments $D_{i, n}^{\ms e}(E)$ and $D_{i, n}^M(E)$

 %
 %MOM BOUND
 %
 \bel[Moment bounds]
 \label{mom_prop}
 Let $k \ge 1$. Then, 
 \been
 \im
$\sup_{n\ge i \ge 1}\sup_{E \su [0, S]}{\E[|D_{i, n}^{\ms e, M}(E)|^k]} < \ff.$
 \im
$\sup_{n\ge i \ge 1}\sup_{E \su \ot}{|E|^{-5/8}}{\E[|D_{i, n}^M(E)|^k]} < \ff.$
 \enen
 \enl

The second essential tool is a rapid decay of covariances between contributions computed from distant regions in space. A similar result appeared in \cite[Lemma 2]{cyl} but was restricted to the assumption that the underlying process is a Poisson point process. In the Gibbs setting, we rely on the idea of exponential clustering appearing in \cite[Lemma 3.4]{gibbs}.

\bel[Covariance bound]
\label{covBoundLem}
There exist constants $C, C' > 0$ with the following property. Let $n \ge 1$ and $S_1, S_2 \su  \on$  with $\#S_1  \vee \#S_2 \le 4$, and set          $X_1 = \prod_{z_i \in S_1}D_{i, n}(E)$ and $X_2 = \prod_{z_j \in S_2 }D_{j, n}(E)$. Then,
$$\Cov\big(X_1, X_2\big) \le C \exp\big(-\dist(S_1, S_2)^{C'}\big)\sqrt{\E[X_1^4]\E[X_2^4]}.$$
\enl

Once Lemmas \ref{mom_prop} and \ref{covBoundLem} are established, the proof of Proposition \ref{var_prop} is obtained by repeating the arguments in \cite[Proposition 3]{cyl}. Instead of doing so, we now elucidate in more detail how the proofs of Lemmas \ref{mom_prop} and \ref{covBoundLem} differ from those encountered in \cite{cyl}.
A key ingredient in the proof of Lemma \ref{mom_prop} is the following variant of \cite[Proposition 2]{cyl}, which can be seen as an extension of Lemma \ref{cont_v_lem} to pairs of filtration times. 
%
%CONT LEM 
%
\bel[Continuity of filtration time pairs]
\label{cont2_lem}
Let $-1 \le q' < q + 1 \le d$ and $\{\eta_ i\}_{i \ge 0}$, $\{\bar\eta_i\}_{i \ge 0}$ be iid uniform in $B_{h_0}$ for some $h_0  > 0$.  Then, there exists a constant $\cv' > 0$ such that for all $P_0, \dots, P_{q + 1}, \bar P_{q' + 1}, \dots, \bar P_q \in \R^p$ and all blocks $E \su [0, S]^2$,
$$\P\big((\rfc(\s),\rfc(\wt \s)) \in E, \rfc(\wt\s) > \rfc(\s)\big) \le \cv'|E|^{3/4},$$
where $\s := \big\{P_0 + \eta_0, \dots, P_{q'} + \eta_{q'}, \bar P_{q' + 1} + \bar\eta_{q' + 1}, \dots, \bar P_q + \bar\eta_q  \big\}$ and  $\wt \s := \big\{P_0 + \eta_ 0, \dots, P_{q + 1} + \eta_{q + 1}\big\}.$
\enl
\bep
As in Lemma \ref{cont_v_lem}, the statement follows from \cite[Proposition 2]{cyl} after expressing the uniform distribution in integral form.
\enp

%NEW_REP
To prove Lemma \ref{mom_prop}, we now provide an alternative representation of the martingale differences $D_{i, n}^{\ms e, M}$. To that end, we let $(\r_{i, n}(t))_{t \in \R}$ be a modification of the original free birth-and-death process $(\r_{n}(t))_{t \in \R}$ such that i) $(\r_{i, n}(t))_{t \in \R}$ coincides with $(\r_{n}(t))_{t \in \R}$ in $W_n \sm W(z_i)$, and ii) $(\r_{i, n}(t))_{t \in \R}$ coincides with an independent copy of $(\r_{W(z_i)}(t))_{t \in \R}$ in $W(z_i)$. Then, we write $\b_{i, n}^{\ms e, M, s}$ and  $\b_{i, n}^{\ms e, M}(E)$ for the persistent Betti number and the increment computed with respect to the Gibbs process obtained from $(\r_{i, n}(t))_{t \in \R}$. We also set $\De_{i, n}^{\ms e, M}(E) := |\b_n^{\ms e, M}(E) - \b_{i, n}^{\ms e, M}(E)|$. In particular, $\E[\b_n^{\ms e, M}(E)\ba \GG_{i - 1}] = \E[\b_{i, n}^{\ms e, M}(E) \ba \GG_i]$  so that $D_{i, n}^{\ms e, M} \le \E[\De_{i, n}^{\ms e, M}(E)\ba \GG_i]$. This notation is adapted accordingly for the $M$-localized  persistent Betti numbers.
%
%MOM PROOF
%
\bep[Proof of Lemma \ref{mom_prop}] We first describe the general strategy, which is common to both parts; hence, we drop the superscripts $\ms e$ and $M$ for the common part. By Jensen's inequality it suffices to establish the moment bounds with $D_{i, n}(E)$ replaced by $\De_{i, n}(E)$. Since $\E[\De_{i, n}(E)^k] = \int_0^\ff k s^{k - 1} \P(\De_{i, n}(E) > s) \d s$, it suffices to show that 
\been
\im[a)] there exists $c > 0$ such that $\P(\De_{i, n}^{}(E) > s) \le  e^{-cs^{c}}$ for all $s \ge 0$; and 
\im[b)] there exists $c' > 0$ such that $\P(\De_{i, n}^{}(E) \ge 1) \le c'|E|^{11/16}$.
\enen 

 {\bf a)} We set $s' := s^{1/(2p(p + 1))}$ and note that $\P(\Rce_n(z_i) > s')$ decays at stretched exponential speed in $s$. Moreover, under the event $\{\Rce_n(z_i) \le s'\}$, the tessellations built from $(\r_n(t))_{t \in \R}$ and from $(\r_{i, n}(t))_{t \in \R}$ coincide in $W_n \sm W_{s'}(z_i)$. Hence, any 2-face that is contained in only one of the tessellations built from $(\r_n(t))_{t \in \R}$ or $(\r_{i, n}(t))_{t \in \R}$ lies within $W_{s'}(z_i)$. Thus,
$$\P\big(\De_{i, n}^{\ms e, M}(E) \ge s, \Rce_n(z_i) \le s'\big) \le \P\big(\Xi^{(2)}_n(W_{s'}(z_i)) \ge s\big),$$
where by Lemma \ref{vert_lem}, the right-hand side decays at stretched exponential speed in $s$. The reasoning for $\De_{i, n}^{M}(E)$ is similar as we essentially need to replace $2$-faces by $q$-faces.

 \medskip

 {\bf b)} By exponential stabilization, it suffices to prove the bound with the  left-hand side replaced by $\P\big(\De_{i, n}(E) \ge 1, \Rce_n(z_i) \le |E|^{-\e_0}\big)$, where $\e_0 := (8(p + 2))^{-3}$.  Next, under the event $\{\Rce_n(z_i) \le |E|^{-\e_0}\}$, the tessellations built from $(\r_n(t))_{t \in \R}$ and from $(\r_{i, n}(t))_{t \in \R}$ coincide in $W_n \sm W_{|E|^{-\e_0}}(z_i)$. We now describe how to proceed for $\De_{i, n}^M(E)$. The arguments for $\De_{i, n}^{\ms e, M}(E)$ are very similar, essentially replacing $q$-simplices by 2-faces.
 Now, write $K_{i, q, n}$ and $K_{i, q, n}'$ for the family of all $(q + 1)$-tuples of points from the vertex sets $\YY_{i, n} := \Xi^{(0)}_n \cap W_{|E|^{-\e_0}}(z_i)$ and $\YY_{i, n}':=\Xi^{(0)}_{i, n} \cap W_{|E|^{-\e_0}}(z_i)$, respectively. Thus, if $\De_{i, n}(E) \ge 1$, then at least one of $K_{i, q, n} \ti K_{i, q + 1, n}$ and $K_{i, q, n}' \ti K_{i, q + 1, n}'$ contains a pair of a $q$-simplex $\sb$ and a $(q + 1)$-simplex $\sd$ such that $\big(\rfc(\sb), \rfc(\sd)\big) \in E$. Therefore, by Lemmas \ref{vert_lem} and \ref{cont2_lem}, and the Cauchy-Schwarz inequality,
	\begin{align*}
		\P\big(\De_{i, n}(E) \ge 1, \Rce_n(z_i) \le |E|^{-\e_0}\big) &\le \E\Big[\sum_{\substack{\sb \in K_{i, q, n} \\ \sd \in K_{i, q + 1, n}}}\P\Big(\big(\rfc(\sb), \rfc(\sd)\big) \in E\ba \XX_n\Big)\Big]\\
		&\le 2\cv'|E|^{3/4}\sqrt{\E\big[|K_{i, q, n}|^2]}\sqrt{\E\big[|K_{i, q + 1, n}|^2]}\\
		&\le 2\cv'|E|^{3/4}\sqrt{\E\big[|\YY_{i, n}|^{2q + 2}]}\sqrt{\E\big[|\YY_{i + 1, n}|^{2q + 4}]}\\
		&\le 2c|E|^{3/4}|E|^{-4\e_0(p + 2)^3},
	\end{align*}
	where we recall that $(\r_n(t))_{t \in \R}$ and $(\r_{i, n}(t))_{t \in \R}$ have the same distribution. Thus, inserting the definition of $\e_0$ concludes the proof.
\enp

%
%COV BOUND PRF
%
\bep[Proof of Lemma \ref{covBoundLem}]
%RE-REP
We first give a different representation of the random variables $X_1, X_2$ that does not rely on computing conditional expectations.
To that end, we first let $\{\r_n^k(\cdot)\}_{k \ge 1}$ of $\r_n$ be independent copies of the free birth-and-death process $\r_n(\cdot)$. Then, we set
        \begin{align*}
		\r_{i, n}^* &:= (\r_n \cap \bigcup_{j \le i}W(z_j)) \cup (\r_n^i \cap \bigcup_{j > i} W(z_j)),\\
		\r_{i, n}^{**} &:= (\r_n \cap \bigcup_{j < i}W(z_j)) \cup (\r_n^i \cap \bigcup_{j \ge i} W(z_j)),\\
		X_\ell' &:= \prod_{i \in S_\ell}\big(\b_n\big(E, \r_{i, n}^*\big) - \b_n\big(E, \r_{i, n}^{**}\big)\big),
	\end{align*}
so that by construction,  $\ms{Cov}(X_1', X_2') = \ms{Cov}(X_1, X_2)$.Now, we let $\Rce_n(S_\ell, \rho_{i, n}^*)$ be the external stabilization radius computed with respect to $\r_{i, n}^*$ and put
$$\Rce_n(S_\ell) := \max_{i \in S_\ell}\big(\Rce(S_\ell, \r_{i, n}^*) \vee \Rce(S_\ell, \r_{i, n}^{**})\big).$$
Set $a := \dist(S_1, S_2)/(8 \sqrt p)$. Then, we decompose $\Cov(X_1, X_2)$ as 
\begin{align*}\Cov(X_1', X_2') &= \Cov(X_1' \one\{\Rce_n(S_1)\le a\}, X_2'\one\{\Rce_n(S_2)\le a\}) \\
	&\phantom= + \Cov(X_1' \one\{\Rce_n(S_1)> a\}, X_2'\one\{\Rce_n(S_2)\le a\})  + \Cov(X_1', X_2'\one\{\Rce_n(S_2)>a\}), 
\end{align*}
where we compute $\Rce_n(S_\ell)$ with respect to $\r_{i_\ell, n}$ for some fixed $i_\ell \in S_\ell$. By the Cauchy-Schwarz inequality, the covariances in the last line can both be bounded by 
$$\E[X_1'X_2' \one\{\Rce_n(S_\ell)> a\}] \le \big(\E[(X_1')^4]\E[(X_2')^4]\big)^{1/4} \sqrt{\P(\Rce_n(S_\ell)> a)},$$
for some $\ell \in\{1, 2\}$. Finally, by the definition of the combined radius of stabilization, we have that $X_\ell' \one \{\Rce_n(S_\ell) \le a\}$, $\ell = 1,2$ are independent so that 
$$\Cov(X_1' \one\{\Rce_n(S_1)\le a\}, X_2'\one\{\Rce_n(S_2)\le a\}) = 0.$$
\enp

\end{document}